%
%

\newif\ifblackboardbold

\blackboardboldtrue
\overfullrule=0pt

\def\bkt{\it} 

\font\titlefont=cmbx12 scaled\magstephalf
\font\sectionfont=cmbx12

\font\foott=cmtt9
\font\footfont=cmr8
\font\footfonts=cmr6

\font\scriptit=cmti10 at 7pt
\font\scriptsl=cmsl10 at 7pt
\scriptfont\itfam=\scriptit
\scriptfont\slfam=\scriptsl


\newfam\bboldfam
\ifblackboardbold
\font\tenbbold=msbm10
\font\sevenbbold=msbm7
\font\fivebbold=msbm5
\textfont\bboldfam=\tenbbold
\scriptfont\bboldfam=\sevenbbold
\scriptscriptfont\bboldfam=\fivebbold
\def\bbold{\fam\bboldfam\tenbbold}
\else
\def\bbold{\bf}
\fi


\newfam\msamfam
\font\tenmsam=msam10
\font\sevenmsam=msam7
\font\fivemsam=msam5
\textfont\msamfam=\tenmsam
\scriptfont\msamfam=\sevenmsam
\scriptscriptfont\msamfam=\fivemsam

\newfam\msbmfam
\font\tenmsbm=msam10
\font\sevenmsbm=msam7
\font\fivemsbm=msam5
\textfont\msbmfam=\tenmsbm
\scriptfont\msbmfam=\sevenmsbm
\scriptscriptfont\msbmfam=\fivemsbm

\newcount\amsfamcount 
\newcount\classcount   
\newcount\positioncount
\newcount\codecount
\newcount\nct             
\def\newsymbol#1#2#3#4#5{               
\nct="#2                                  
\ifnum\nct=1 \amsfamcount=\msamfam\else   
\ifnum\nct=2 \amsfamcount=\msbmfam\else   
\ifnum\nct=3 \amsfamcount=\eufmfam
\fi\fi\fi
\multiply\amsfamcount by "100           
\classcount="#3                 
\multiply\classcount by "1000           
\positioncount="#4#5            
\codecount=\classcount                  
\advance\codecount by \amsfamcount      
\advance\codecount by \positioncount
\mathchardef#1=\codecount}              


\font\Arm=cmr9
\font\Ai=cmmi9
\font\Asy=cmsy9
\font\Abf=cmbx9
\font\Brm=cmr7
\font\Bi=cmmi7
\font\Bsy=cmsy7
\font\Bbf=cmbx7
\font\Crm=cmr6
\font\Ci=cmmi6
\font\Csy=cmsy6
\font\Cbf=cmbx6

\ifblackboardbold
\font\Abbold=msbm10 at 9pt
\font\Bbbold=msbm7
\font\Cbbold=msbm5 at 6pt
\fi


\font\small=cmr5

\relax



\def\small{%
\textfont0=\Arm \scriptfont0=\Brm \scriptscriptfont0=\Crm
\textfont1=\Ai \scriptfont1=\Bi \scriptscriptfont1=\Ci
\textfont2=\Asy \scriptfont2=\Bsy \scriptscriptfont2=\Csy
\textfont\bffam=\Abf \scriptfont\bffam=\Bbf \scriptscriptfont\bffam=\Cbf
\def\rm{\fam0\Arm}\def\mit{\fam1}\def\oldstyle{\fam1\Ai}%
\def\bf{\fam\bffam\Abf}%
\ifblackboardbold
\textfont\bboldfam=\Abbold
\scriptfont\bboldfam=\Bbbold
\scriptscriptfont\bboldfam=\Cbbold
\def\bbold{\fam\bboldfam\Abbold}%
\fi
\rm
}

\def\ifootnote#1#2{\footnote{#1}{{\lipmanpoint\rm #2\vskip-4ex\egroup}}}

\font\tenmsb=msbm10
\font\sevenmsb=msbm10 at 7pt
\font\fivemsb=msbm10 at 5pt
\newfam\msbfam
\textfont\msbfam=\tenmsb
\scriptfont\msbfam=\sevenmsb
\scriptscriptfont\msbfam=\fivemsb

\def\hexnumber#1{\ifcase#1 0\or1\or2\or3\or4\or5\or6\or7\or8\or9\or
	A\or B\or C\or D\or E\or F\fi}

\mathchardef\subsetneq="2\hexnumber\msbfam28


\ifx\begin\undefined\else\global\advance\srcdepth by
1\expandafter\endinput\fi

\def\begin{}
\newcount\srcdepth
\srcdepth=1

\outer\def\bye{\global\advance\srcdepth by -1
  \ifnum\srcdepth=0
    \def\endcmd{\vfill\supereject\nopagenumbers\par\vfill\supereject\end}
  \else\def\endcmd{}\fi
  \endcmd
}



\def\initialize#1#2#3#4#5#6{
  \ifnum\srcdepth=1
  \magnification=#1
  \hsize = #2
  \vsize = #3
  \hoffset=#4
  \advance\hoffset by -\hsize
  \divide\hoffset by 2
  \advance\hoffset by -1truein
  \voffset=#5
  \advance\voffset by -\vsize
  \divide\voffset by 2
  \advance\voffset by -1truein
  \advance\voffset by #6
  \baselineskip=13pt
  \emergencystretch = 0.05\hsize
  \fi
}

\def\irenasize{\initialize{1195}
  {6.5truein}{8.8truein}{8.5truein}{11truein}{0.2truein}}








\def\quiremacro#1#2#3#4#5#6#7#8#9{
  \expandafter\def\csname#1\endcsname##1{
  \ifnum\srcdepth=1
  \magnification=#2
  \input quire
  \hsize=#3
  \vsize=#4
  \htotal=#5
  \vtotal=#6
  \shstaplewidth=#7
  \shstaplelength=#8
  \hoffset=\htotal
  \advance\hoffset by -\hsize
  \divide\hoffset by 2
  \ifnum\vsize<\vtotal
    \voffset=\vtotal
    \advance\voffset by -\vsize
    \divide\voffset by 2
  \fi
  \advance\voffset by #9
  \shhtotal=2\htotal
  \baselineskip=13pt
  \emergencystretch = 0.05\hsize
  \horigin=0.0truein
  \vorigin=0.0truein
  \shthickness=0pt
  \shoutline=0pt
  \shcrop=0pt
  \shvoffset=-1.0truein
  \ifnum##1>0\quire{#1}\else\qtwopages\fi
  \fi
}}



\quiremacro{letterbooklet} 
{1000}{4.79452truein}{7truein}{5.5truein}{8.5truein}{0.01pt}{0.66truein}
{-.0625truein}

\quiremacro{Afourbooklet}
{1095}{5.25truein}{7truein}{421truept}{595truept}{0.01pt}{0.66truein}
{-.0625truein}

\quiremacro{legalbooklet}
{1095}{5.25truein}{7truein}{7.0truein}{8.5truein}{0.01pt}{0.66truein}
{-.0625truein}

\quiremacro{twoupsub} 
{895}{4.5truein}{7truein}{5.5truein}{8.5truein}{0pt}{0pt}{.0625truein}


\quiremacro{Afourviewsub} 
{1000}{5.0228311in}{7.7625571in}{421truept}{595truept}{0.1pt}{0.5\vtotal}
{-.0625truein}


\quiremacro{viewsub}
{1095}{5.5truein}{8.5truein}{461truept}{666truept}{0.1pt}{0.5\vtotal}
{-.125truein}



%






\newlinechar=`@
\def\forwardmsg#1#2#3{\immediate\write16{@*!*!*!* forward reference should
be: @\noexpand\forward{#1}{#2}{#3}@}}
\def\nodefmsg#1{\immediate\write16{@*!*!*!* #1 is an undefined reference@}}

\def\forwardsub#1#2{\def\newref{{#2}{#1}}}

\def\forward#1#2#3{%
\expandafter\expandafter\expandafter\forwardsub\expandafter{#3}{#2}
\expandafter\ifx\csname#1\endcsname\relax\else%
\expandafter\ifx\csname#1\endcsname\newref\else%
\expandafter\ifx\csname#2\endcsname\relax\else%
\forwardmsg{#1}{#2}{#3}\fi\fi\fi%
\expandafter\let\csname#1\endcsname\newref}

\def\firstarg#1{\expandafter\argone #1}\def\argone#1#2{#1}
\def\secondarg#1{\expandafter\argtwo #1}\def\argtwo#1#2{#2}

\def\ref#1{\expandafter\ifx\csname#1\endcsname\relax
  {\nodefmsg{#1}\bf`#1'}\else
  \expandafter\firstarg\csname#1\endcsname
  ~\htmllocref{#1}{\expandafter\secondarg\csname#1\endcsname}\fi}

\def\refscor#1{\expandafter\ifx\csname#1\endcsname\relax
  {\nodefmsg{#1}\bf`#1'}\else
  Corollaries~\htmllocref{#1}{\expandafter\secondarg\csname#1\endcsname}\fi}

\def\refs#1{\expandafter\ifx\csname#1\endcsname\relax
  {\nodefmsg{#1}\bf`#1'}\else
  \expandafter\firstarg\csname #1\endcsname
  s~\htmllocref{#1}{\expandafter\secondarg\csname#1\endcsname}\fi}

\def\refn#1{\expandafter\ifx\csname#1\endcsname\relax
  {\nodefmsg{#1}\bf`#1'}\else
  \htmllocref{#1}{\expandafter\secondarg\csname #1\endcsname}\fi}

\def\pageref#1{\expandafter\ifx\csname#1\endcsname\relax%
  {\nodefmsg{#1}\bf`#1'}\else%
  \htmllocref{#1}{\expandafter\csname#1\endcsname}\fi}

\def\slidepageref#1{\expandafter\ifx\csname#1\endcsname\relax%
  {\nodefmsg{#1}\bf`#1'}\else%
  \htmllocref{#1}{\expandafter\csname#1\endcsname}{1}\fi}

\def\slidepageref#1{\expandafter\ifx\csname#1\endcsname\relax%
  {\nodefmsg{#1}\bf`#1'}\else%
  \htmllocref{#1}{\expandafter\csname#1\endcsname}{ }\fi}

\def\blueslidepageref#1{\expandafter\ifx\csname#1\endcsname\relax%
  {\nodefmsg{#1}\bf`#1'}\else%
  \htmllocref{#1}{\expandafter\csname#1\endcsname}{2}\fi}

\def\blueslidepageref#1{\expandafter\ifx\csname#1\endcsname\relax%
  {\nodefmsg{#1}\bf`#1'}\else%
  \htmllocref{#1}{\expandafter\csname#1\endcsname}{ }\fi}

\def\slidepagelabel#1#2{%
   \write\isauxout{\noexpand\forward{\noexpand#1}{}{#2}}%
   \bgroup\htmlanchor{#1}{\makeref{\noexpand#1}{}{{#2}}}\egroup%
 }






%
%

%
%
%
%
%
%
%
%


\newif\ifhypers  
\hypersfalse


\edef\freehash{\catcode`\noexpand\#=\the\catcode`\#}%
\catcode`\#=12
\freehash
\let\freehash=\relax
\ifhypers\fi
\def\puthtml#1{\ifhypers\fi}
\def\htmlanchor#1#2{\puthtml{<a name="#1">}#2\puthtml{</a>}}

\def\@pdfm@mark#1{}
\def\setlink#1{\colored{\linkcolor}{#1}}%
\def\setlink#1{\Purple{#1}}%

%

%
\def\htmllocref#1#2{\ifhypers\leavevmode\fi\setlink{#2}\ifhypers\fi\relax}%
%
%
%
%
\def\Acrobatmenu#1#2{%
  \@pdfm@mark{%
    bann <<
      /Type /Annot
      /Subtype /Link
      /A <<
        /S /Named
        /N /#1
      >>
      /Border [\@pdfborder]
      /C [\@menubordercolor]
    >>%
   }%
  \Hy@colorlink{\@menucolor}#2\Hy@endcolorlink
  \@pdfm@mark{eann}%
}
\def\@pdfborder{0 0 1}
\def\@menubordercolor{1 0 0}
\def\@menucolor{red}

\def\ifempty#1#2\endB{\ifx#1\endA}
\def\makeref#1#2#3{\ifempty#1\endA\endB\else\forward{#1}{#2}{#3}\fi\unskip}


\def\crosshairs#1#2{
  \dimen1=.002\drawx
  \dimen2=.002\drawy
  \ifdim\dimen1<\dimen2\dimen3\dimen1\else\dimen3\dimen2\fi
  \setbox1=\vclap{\vline{2\dimen3}}
  \setbox2=\clap{\hline{2\dimen3}}
  \setbox3=\hstutter{\kern\dimen1\box1}{4}
  \setbox4=\vstutter{\kern\dimen2\box2}{4}
  \setbox1=\vclap{\vline{4\dimen3}}
  \setbox2=\clap{\hline{4\dimen3}}
  \setbox5=\clap{\copy1\hstutter{\box3\kern\dimen1\box1}{6}}
  \setbox6=\vclap{\copy2\vstutter{\box4\kern\dimen2\box2}{6}}
  \setbox1=\vbox{\offinterlineskip\box5\box6}
  \smash{\vbox to #2{\hbox to #1{\hss\box1}\vss}}}

\def\boxgrid#1{\rlap{\vbox{\offinterlineskip
  \setbox0=\hline{\wd#1}
  \setbox1=\vline{\ht#1}
  \smash{\vbox to \ht#1{\offinterlineskip\copy0\vfil\box0}}
  \smash{\vbox{\hbox to \wd#1{\copy1\hfil\box1}}}}}}

\def\drawgrid#1#2{\vbox{\offinterlineskip
  \dimen0=\drawx
  \dimen1=\drawy
  \divide\dimen0 by 10
  \divide\dimen1 by 10
  \setbox0=\hline\drawx
  \setbox1=\vline\drawy
  \smash{\vbox{\offinterlineskip
    \copy0\vstutter{\kern\dimen1\box0}{10}}}
  \smash{\hbox{\copy1\hstutter{\kern\dimen0\box1}{10}}}}}

\long\def\boxit#1#2#3{\hbox{\vrule
  \vtop{%
    \vbox{\hsize=#2\hrule\kern#1%
      \hbox{\kern#1#3\kern#1}}%
    \kern#1\hrule}%
    \vrule}}

\newdimen\boxrulethickness \boxrulethickness=.4pt
\newdimen\boxedhsize
\newbox\textbox 
\newdimen\originalbaseline 
\newdimen\hborderspace\newdimen\vborderspace
\hborderspace=3pt \vborderspace=3pt

\def\preparerulebox#1#2{\setbox\textbox=\hbox{#2}%
   \originalbaseline=\vborderspace
   \advance\originalbaseline\boxrulethickness
   \advance\originalbaseline\dp\textbox
   \def\Borderbox{\vbox{\hrule height\boxrulethickness
     \hbox{\vrule width\boxrulethickness\hskip\hborderspace
     \vbox{\vskip\vborderspace\relax#1{#2}\vskip\vborderspace}%
     \hskip\hborderspace\vrule width\boxrulethickness}%
     \hrule height\boxrulethickness}}}

\def\hrulebox#1{\hbox{\preparerulebox{\hbox}{#1}%
   \lower\originalbaseline\Borderbox}}
\def\vrulebox#1#2{\vbox{\preparerulebox{\vbox}{\hsize#1#2}\Borderbox}}
\def\parrulebox#1\par{\boxedhsize=\hsize\advance\boxedhsize
   -2\boxrulethickness \advance\boxedhsize-2\hborderspace
   \vskip3\parskip\vrulebox{\boxedhsize}{#1}\vskip\parskip\par}
\def\sparrulebox#1\par{\vskip1truecm\boxedhsize=\hsize\advance\boxedhsize
   -2\boxrulethickness \advance\boxedhsize-2\hborderspace
   \vskip3\parskip\vrulebox{\boxedhsize}{#1}\vskip.5truecm\par}



%

\newcount\sectno \sectno=0
\newcount\appno \appno=64 
\newcount\thmno \thmno=0
\newcount\randomct \newcount\rrandomct
\newif\ifsect \secttrue
\newif\ifapp \appfalse


\def\widow#1{\vskip 0pt plus#1\vsize\goodbreak\vskip 0pt plus-#1\vsize}


\def\stdskip{\vskip 9pt plus3pt minus 3pt}
\def\stdbreak{\par\removelastskip\penalty-100\stdskip}

\def\proof{\stdbreak\noindent \colored{\proofcolor}{{\sl Proof. }}}
\def\proof{\stdbreak\noindent \ForestGreen{{\sl Proof. }}}


\newif\ifnumberbibl \numberbibltrue
\newdimen\labelwidth

\def\references{\bgroup
  \edef\numtoks{}%
  \global\thmno=0
  \setbox1=\hbox{[999]} 
  \labelwidth=\wd1 \advance\labelwidth by 2.5em
  \ifnumberbibl\advance\labelwidth by -2em\fi
  \parindent=\labelwidth \advance\parindent by 0.5em 
  \removelastskip
  \widow{0.1}
  \vskip 24pt plus 6pt minus 6 pt
  \frenchspacing
  \immediate\write\isauxout{\noexpand\forward{Bibliography}{}{\the\pageno}}%
  \immediate\write\iscontout{\noexpand\contnosectlist{}{Bibliography}{\the\pageno}}%
  \leftline{\sectionfont \Orange{References}}
  \ifhypers%
     \global\thmno=0\relax%
     \ifsect%
	\global\advance\sectno by 1%
	\edef\numtoks{\number\sectno}%
  	\hbox{}%
     \else%
	\edef\numtoks{}%
  	\hbox{}%
     \fi%
  \fi%
  \nobreak\stdskip\noindent}%
\def\endreferences{\nonfrenchspacing\egroup}

\def\bitem#1{\global\advance\thmno by 1%
  \outer\expandafter\def\csname#1\endcsname{\the\thmno}%
  \edef\numtoks{\number\thmno}%
  \ifhypers\htmlanchor{#1}{\makeref{\noexpand#1}{REF}{\numtoks}}\fi%
  \ifnumberbibl
    \immediate\write\isauxout{\noexpand\forward{\noexpand#1}{}{\the\thmno}}
    \item{\hbox to \labelwidth{\Maroon{\hfil\the\thmno.\ \ }}}
  \else
    \immediate\write\isauxout{\noexpand\expandafter\noexpand\gdef\noexpand\csname #1\noexpand\endcsname{#1}}
    \item{\hbox to \labelwidth{\Olive{\hfil#1\ \ }}}
  \fi}

\outer\def\section#1#2{%
  \appfalse\secttrue%
  \removelastskip%
  \global\advance\sectno by 1%
  \global\thmno=0\relax%
  \edef\numtoks{\number\sectno}%
  \vskip 24pt plus 6pt minus 6 pt%
  \widow{.03}%
  \leftline{{\draftlabel{#1}\sectionfont\OliveGreen{\numtoks}\quad \Apricot{#2}}}%
  \immediate\write\isauxout{\noexpand\forward{\noexpand#1}{Section}{\numtoks}}%
  \immediate\write\iscontout{\noexpand\contlist{\noexpand#1}{#2}{\the\pageno}}%
  \htmlanchor{#1}{\makeref{#1}{Section}{\numtoks}}%
  \message{#2}%
  \ifhypers\hbox{}\fi%
  \nobreak\vskip 3pt}

\outer\def\appendix#1#2{%
  \removelastskip%
  \global\advance\appno by 1%
  \global\thmno=0\relax%
  \edef\numtoks{Appendix\ \char\number\appno}%
  \vskip 24pt plus 6pt minus 6 pt%
  \widow{.03}%
  \leftline{\draftlabel{#1}\sectionfont\NavyBlue{\numtoks}\quad \Blue{#2}}%
  \immediate\write\isauxout{\noexpand\forward{\noexpand#1}{}{\numtoks}}%
  \immediate\write\iscontout{\noexpand\contnosectlist{}{\numtoks:  #2}{\the\pageno}}%
  \htmlanchor{#1}{\makeref{#1}{Section}{\numtoks}}%
  \message{#2}%
  \ifhypers\hbox{}\fi%
  \apptrue\sectfalse%
  \nobreak\vskip 3pt}


\def\proclamation#1#2#3{
  \outer\expandafter\def\csname#1\endcsname{\bgroup%
    \gdef\Prnm{#2}%
    \global\advance\thmno by 1%
    \global\randomct=0%
    \ifsect\gdef\numtoks{\number\sectno.\number\thmno}%
    \else\ifapp\gdef\numtoks{\char\number\appno.\number\thmno}%
    \else\gdef\numtoks{\number\thmno #3}%
    \fi\fi%
    \stdbreak%
    \noindent{\bf \Fuchsia{#2\ \numtoks}\enspace}%
    \futurelet\testchar\maybeoptionproclaim}}

\def\maybeoptionproclaim{\ifx[\testchar\let\next\optionproclaim%
	\else\let\next\nooptionproclaim\fi\next}

\def\optionproclaim[#1]{%
  \htmlanchor{#1}{\makeref{\noexpand#1}{\Prnm}{\numtoks}}%
  \immediate\write\isauxout{\noexpand\forward{\noexpand#1}{\Prnm}{\numtoks}}%
  \draftlabel{#1}%
  \sl}
\def\nooptionproclaim{\sl}
\def\endb{\par\stdbreak\egroup\randomct=0}

\def\pagelabel#1{%
   \write\isauxout{\noexpand\forward{\noexpand#1}{page}{\ \the\pageno}}%
   \bgroup\htmlanchor{#1}{\makeref{\noexpand#1}{page}{\ \the\pageno}}\egroup%
   \draftlabel{page\ #1}%
 }

\def\eqlabel#1{
    \global\advance\thmno by 1%
    \ifsect\gdef\numtoks{(\number\sectno.\number\thmno)}%
    \else\ifapp\gdef\numtoks{(\char\number\appno.\number\thmno)}%
    \else\gdef\numtoks{(\number\thmno)}%
    \fi\fi%
    \outer\expandafter\def\csname#1\endcsname{\numtoks}%
    \htmlanchor{#1}{\makeref{\noexpand#1}{}{\numtoks}}%
    \immediate\write\isauxout{\noexpand\forward{\noexpand#1}{equation}{\numtoks}}%
    \quad\qquad\hfill\eqno\hbox{\hfill\rm\noexpand\Green{\numtoks}\draftcmt{\OliveGreen{#1}}}%
}
\def\eqalabel#1{
    \global\advance\thmno by 1%
    \ifsect\gdef\numtoks{(\number\sectno.\number\thmno)}%
    \else\ifapp\gdef\numtoks{(\char\number\appno.\number\thmno)}%
    \else\gdef\numtoks{(\number\thmno)}%
    \fi\fi%
    \outer\expandafter\def\csname#1\endcsname{\numtoks}%
    \htmlanchor{#1}{\makeref{\noexpand#1}{}{\numtoks}}%
    \immediate\write\isauxout{\noexpand\forward{\noexpand#1}{equation}{\numtoks}}%
    \quad\qquad&\hfill\hbox{\hfill\rm\noexpand\Green{\numtoks}\draftcmt{\OliveGreen{#1}}}%
}


\def\othernumbered#1#2{
  \outer\expandafter\def\csname#1\endcsname{%
  \gdef\Prnm{#2}%
  \global\advance\thmno by 1%
  \global\randomct=0%
  \ifsect\gdef\numtoks{\number\sectno.\number\thmno}%
  \else\gdef\numtoks{\number\thmno}%
  \fi%
  \stdbreak%
  \noindent{\bf \Bittersweet{#2\ \numtoks}\enspace}%
  \futurelet\testchar\maybeoptionothernumbered}}

\def\stothernumbered#1#2{
  \outer\expandafter\def\csname#1\endcsname{%
  \gdef\Prnm{#2}%
  \global\advance\thmno by 1%
  \global\randomct=0%
  \ifsect\gdef\numtoks{\number\sectno.\number\thmno}%
  \else\gdef\numtoks{\number\thmno}%
  \fi%
  \stdbreak%
  \noindent{\bf \Bittersweet{#2\ \numtoks*}\enspace}%
  \futurelet\testchar\maybeoptionothernumbered}}

\def\impnumbered#1#2{
  \outer\expandafter\def\csname#1\endcsname{%
  \gdef\Prnm{#2}%
  \global\advance\thmno by 1%
  \global\randomct=0%
  \ifsect\gdef\numtoks{\number\sectno.\number\thmno}%
  \else\gdef\numtoks{\number\thmno}%
  \fi%
  \stdbreak%
  \noindent{\bf \Bittersweet{#2\ \numtoks\ (IMPORTANT!)}\enspace}%
  \futurelet\testchar\maybeoptionothernumbered}}

\def\maybeoptionothernumbered{\ifx[\testchar\let\next\optionothernumbered%
	\else\let\next\nooptionothernumbered\fi\next}

\def\optionothernumbered[#1]{{%
  \htmlanchor{#1}{\makeref{\noexpand#1}{\Prnm}{\numtoks}}%
  \immediate\write\isauxout{\noexpand\forward{\noexpand#1}{\Prnm}{\numtoks}}%
  \draftlabel{#1}}}
\def\nooptionothernumbered{}

\def\label#1{\optionothernumbered{#1}}

\proclamation{definition}{Definition}{}
\proclamation{defin}{Definition}{}
\proclamation{lemma}{Lemma}{}
\proclamation{proposition}{Proposition}{}
\proclamation{prop}{Proposition}{}
\proclamation{theorem}{Theorem}{}
\proclamation{thm}{Theorem}{}
\proclamation{corollary}{Corollary}{}
\proclamation{cor}{Corollary}{}
\proclamation{conjecture}{Conjecture}{}
\proclamation{proc}{Procedure}{}
\proclamation{assumption}{Assumption}{}
\proclamation{axiom}{Axiom}{}

\othernumbered{example}{Example}
\othernumbered{examples}{Examples}
\othernumbered{remark}{Remark}
\othernumbered{remarks}{Remarks}
\othernumbered{remarksexamples}{Remarks and examples}
\othernumbered{facts}{Facts}
\othernumbered{fact}{Fact}
\othernumbered{irrelevant}{Irrelevant}
\othernumbered{question}{Question}
\othernumbered{exercise}{Exercise}
\othernumbered{importantexercise}{IMPORTANT EXERCISE}
\othernumbered{construction}{Construction}
\othernumbered{problem}{Problem}
\othernumbered{table}{Table}
\othernumbered{notation}{Notation}
\othernumbered{figure}{Figure}

\stothernumbered{stexercise}{Exercise}
\stothernumbered{exercisest}{Exercise}
\impnumbered{importantexercise}{Exercise}

\proclamation{lemmap}{Lemma}{'}
\proclamation{lemmapp}{Lemma}{''}

\def\crosshairs#1#2{
  \dimen1=.002\drawx
  \dimen2=.002\drawy
  \ifdim\dimen1<\dimen2\dimen3\dimen1\else\dimen3\dimen2\fi
  \setbox1=\vclap{\vline{2\dimen3}}
  \setbox2=\clap{\hline{2\dimen3}}
  \setbox3=\hstutter{\kern\dimen1\box1}{4}
  \setbox4=\vstutter{\kern\dimen2\box2}{4}
  \setbox1=\vclap{\vline{4\dimen3}}
  \setbox2=\clap{\hline{4\dimen3}}
  \setbox5=\clap{\copy1\hstutter{\box3\kern\dimen1\box1}{6}}
  \setbox6=\vclap{\copy2\vstutter{\box4\kern\dimen2\box2}{6}}
  \setbox1=\vbox{\offinterlineskip\box5\box6}
  \smash{\vbox to #2{\hbox to #1{\hss\box1}\vss}}}

\def\boxgrid#1{\rlap{\vbox{\offinterlineskip
  \setbox0=\hline{\wd#1}
  \setbox1=\vline{\ht#1}
  \smash{\vbox to \ht#1{\offinterlineskip\copy0\vfil\box0}}
  \smash{\vbox{\hbox to \wd#1{\copy1\hfil\box1}}}}}}

\def\drawgrid#1#2{\vbox{\offinterlineskip
  \dimen0=\drawx
  \dimen1=\drawy
  \divide\dimen0 by 10
  \divide\dimen1 by 10
  \setbox0=\hline\drawx
  \setbox1=\vline\drawy
  \smash{\vbox{\offinterlineskip
    \copy0\vstutter{\kern\dimen1\box0}{10}}}
  \smash{\hbox{\copy1\hstutter{\kern\dimen0\box1}{10}}}}}


\def\hpad#1#2#3{\hbox{\kern #1\hbox{#3}\kern #2}}
\def\vpad#1#2#3{\setbox0=\hbox{#3}\vbox{\kern #1\box0\kern #2}}




\def\stack#1#2#3{\vbox{\offinterlineskip
  \setbox2=\hbox{#2}
  \setbox3=\hbox{#3}
  \dimen0=\ifdim\wd2>\wd3\wd2\else\wd3\fi
  \hbox to \dimen0{\hss\box2\hss}
  \kern #1
  \hbox to \dimen0{\hss\box3\hss}}}


\def\hexp#1{%
  \setbox0=\hbox{${}^{#1}$}%
  \hbox to .5\wd0{\box0\hss}}

\def\hsub#1{%
  \setbox0=\hbox{${}_{#1}$}%
  \hbox to .5\wd0{\box0\hss}}



\def\bmatrix#1#2{{\left(\vcenter{\halign
  {&\kern#1\hfil$##\mathstrut$\kern#1\cr#2}}\right)}}

\def\rightarrowmat#1#2#3{
  \setbox1=\hbox{\small\kern#2$\bmatrix{#1}{#3}$\kern#2}
  \,\vbox{\offinterlineskip\hbox to\wd1{\hfil\copy1\hfil}
    \kern 3pt\hbox to\wd1{\rightarrowfill}}\,}

\def\leftarrowmat#1#2#3{
  \setbox1=\hbox{\small\kern#2$\bmatrix{#1}{#3}$\kern#2}
  \,\vbox{\offinterlineskip\hbox to\wd1{\hfil\copy1\hfil}
    \kern 3pt\hbox to\wd1{\leftarrowfill}}\,}

\def\rightarrowbox#1#2{
  \setbox1=\hbox{\kern#1\hbox{\small #2}\kern#1}
  \,\vbox{\offinterlineskip\hbox to\wd1{\hfil\copy1\hfil}
    \kern 3pt\hbox to\wd1{\rightarrowfill}}\,}

\def\leftarrowbox#1#2{
  \setbox1=\hbox{\kern#1\hbox{\small #2}\kern#1}
  \,\vbox{\offinterlineskip\hbox to\wd1{\hfil\copy1\hfil}
    \kern 3pt\hbox to\wd1{\leftarrowfill}}\,}


\def\mOth{\mathsurround=0pt}
\newdimen\pOrenwd \setbox0=\hbox{\tenex B} \pOrenwd=\wd0
\def\lefttopextramatrix#1{\begingroup \mOth
  \setbox0=\vbox{\def\cr{\crcr\noalign{\kern2pt\global\let\cr=\endline}}
    \ialign{$##$\hfil\kern2pt\kern\pOrenwd&\thinspace\hfil$##$\hfil
      &&\quad\hfil$##$\hfil\crcr
      \omit\strut\hfil\crcr\noalign{\kern-\baselineskip}
      #1\crcr\omit\strut\cr}}
  \setbox2=\vbox{\unvcopy0 \global\setbox1=\lastbox}
  \setbox2=\hbox{\unhbox1 \unskip \global\setbox1=\lastbox}
  \setbox2=\hbox{$\kern\wd1\kern-\pOrenwd \left[ \kern-\wd1
    \global\setbox1=\vbox{\box1\kern2pt}
    \vcenter{\kern-\ht1 \unvbox0 \kern-\baselineskip} \,\right]$}
  \null\;\vbox{\kern\ht1\box2}\endgroup}

\def\mOth{\mathsurround=0pt}
\newdimen\pOrenwd \setbox0=\hbox{\tenex B} \pOrenwd=\wd0
\def\leftextramatrix#1{\begingroup \mOth
  \setbox0=\vbox{\def\cr{\crcr\noalign{\kern2pt\global\let\cr=\endline}}
    \ialign{$##$\hfil\kern2pt\kern\pOrenwd&\thinspace\hfil$##$\hfil
      &&\quad\hfil$##$\hfil\crcr
      \omit\strut\hfil\crcr\noalign{\kern-\baselineskip}
      #1\crcr\omit\strut\cr}}
  \setbox2=\vbox{\unvcopy0 \global\setbox1=\lastbox}
  \setbox2=\hbox{\unhbox1 \unskip \global\setbox1=\lastbox}
  \setbox2=\hbox{$\kern\wd1\kern-\pOrenwd \left[ \kern-\wd1
    \vcenter{\unvbox0 \kern-\baselineskip} \,\right]$}
  \null\;\vbox{\box2}\endgroup}


\newcount\countA
\newcount\countB
\newcount\countC

\def\monthname{\begingroup
  \ifcase\number\month
    \or January\or February\or March\or April\or May\or June\or
    July\or August\or September\or October\or November\or December\fi
\endgroup}

\def\dayname{\begingroup
  \countA=\number\day
  \countB=\number\year
  \advance\countA by 0 
  \advance\countA by \ifcase\month\or
    0\or 31\or 59\or 90\or 120\or 151\or
    181\or 212\or 243\or 273\or 304\or 334\fi
  \advance\countB by -1995
  \multiply\countB by 365
  \advance\countA by \countB
  \countB=\countA
  \divide\countB by 7
  \multiply\countB by 7
  \advance\countA by -\countB
  \advance\countA by 1
  \ifcase\countA\or Sunday\or Monday\or Tuesday\or Wednesday\or
    Thursday\or Friday\or Saturday\fi
\endgroup}

\def\timename{\begingroup
   \countA = \time
   \divide\countA by 60
   \countB = \countA
   \countC = \time
   \multiply\countA by 60
   \advance\countC by -\countA
   \ifnum\countC<10\toks1={0}\else\toks1={}\fi
   \ifnum\countB<12 \toks0={\sevenrm AM}
     \else\toks0={\sevenrm PM}\advance\countB by -12\fi
   \relax\ifnum\countB=0\countB=12\fi
   \hbox{\the\countB:\the\toks1 \the\countC \thinspace \the\toks0}
\endgroup}

\def\timestamp{\dayname, \the\day\ \monthname\ \the\year, \timename}


\def\COMMENT#1\par{\bigskip\hrule\smallskip#1\smallskip\hrule\bigskip}

{\catcode`\^^M=12 \endlinechar=-1 %
 \gdef\xcomment#1^^M{\def\test{#1}
   \ifx\test\endcomment \let\next=\endgroup
   \else \let\next=\xcomment \fi
   \next}
}
\def\dospecials{\do\ \do\\\do\{\do\}\do\$\do\&%
  \do\#\do\^\do\^^K\do\_\do\^^A\do\%\do\~}
\def\makeinnocent#1{\catcode`#1=12}
\def\comment{\begingroup
  \let\do=\makeinnocent \dospecials
  \endlinechar`\^^M \catcode`\^^M=12 \xcomment}
{\escapechar=-1
 \xdef\endcomment{\string\\endcomment}
}


\newif\ifdraft\draftfalse


\def\draftcmt#1{\llap{\smash{\raise2ex\hbox{{#1}}}}}
\def\draftind#1{\llap{\smash{\raise-1.3ex\hbox{{#1}}}}}
\def\draftlabel#1{\ifdraft\draftcmt{\footfonts\Red{{\string#1}\ }\unskip\fi}}
\def\draftlabel#1{%
  \ifx#1\empty\else
  \ifdraft\draftcmt{\foott\Red{\string #1}\ \ }\unskip\fi\fi}

\newcount\hours
\newcount\minutes
\hours\time \divide\hours 60
\minutes-\hours \multiply\minutes 60\advance\minutes\time
\edef\timehhmm{\ifnum\hours<10 0\fi\the\hours
:\ifnum\minutes<10 0\fi\the\minutes}

\footline={{\rm\hfill\folio\hfill}}


\newif\ifcolor\colorfalse
\ifcolor\input colordvi\else\input blackdvi\fi%
\def\setcolor#1{\ifnum#1=1\colortrue \input colordvi
	\else\colorfalse\input blackdvi\fi}

\def\cite{\futurelet\testchar\maybeoptioncite}
\def\maybeoptioncite{\ifx[\testchar \let\next\optioncite
	\else \let\next\nooptioncite \fi \next}

\def\nooptioncite#1{\expandafter\ifx\csname#1\endcsname\relax
  {\nodefmsg{#1}\bf`#1'}\else
  \htmllocref{#1}{[\csname#1\endcsname]}{}\fi}

\def\optioncite[#1]#2{\ifx\csname#1\endcsname\relax\else
  ~\htmllocref{#2}{[\csname#2\endcsname, #1]}{}\fi}

\def\strutdepth{\dp\strutbox}
\def\strutdepth{2ex}
\def\margin#1{\strut\vadjust{\kern-\strutdepth%
        \vtop to\strutdepth{\baselineskip\strutdepth%
        \llap{\vbox{\hsize=3em\noindent #1}\hskip1em}\null}}}
\def\rmargin#1{\strut\vadjust{\kern-\strutdepth%
        \vtop to\strutdepth{\baselineskip\strutdepth%
        \vss\rlap{\ \hskip\hsize\ \vtop{\hsize=4em\noindent #1}}\null}}}

\def\quote{\bgroup\footfont\parindent=0pt%
	\baselineskip=8pt\hangindent=15em\hangafter=0}
\def\endquote{\egroup}

\def\spm{\raise0.5pt\hbox{\bgroup\textfont2=\fivesy\relax$\pm$\egroup}}



\def\buildrelu#1\over#2{\mathrel{\mathop{\kern0pt #1}\limits_{#2}}}

\newdimen\hngind \hngind=\parindent \advance\hngind by 1.9em\relax
\def\nonuritem{{\vskip0pt\hangafter=0\global\hangindent\hngind%
        {\vphantom{(\the\randomct)\hskip.5em}}}}
\def\nonuritem{{\parindent=0pt\hngind=\parindent \advance\hngind by 1.9em%
        \vskip0pt\hangafter=0\global\hangindent\hngind%
        \vphantom{{(\the\randomct)\hskip.5em}}}}
\def\nonuritem{{\idoitem{}}}
\def\rlabel#1{{\global\edef#1{\the\randomct}}\ignorespaces}
\def\pitem{{\global\advance\randomct by 1%
	\vskip1pt\hangafter=0\global\hangindent\hngind%
	{{\the\randomct.\hskip.5em}}}}
\def\rritem{{\global\advance\rrandomct by 1%
	\vskip1pt\hangafter=0\hangindent3\hngind%
	\ \hskip0em{\hbox to2.5em{\hfill(\romannumeral\rrandomct)\hfill}}}}
\def\rlabel#1{{%
   \bgroup\htmlanchor{#1}{\makeref{\noexpand#1}{item}{\ \the\randomct}}\egroup%
   \immediate\write\isauxout{\noexpand\forward{\noexpand#1}{item}{\ \the\randomct}}%
   \draftlabel{item\ #1}}}

\def\idoitem#1{{\parindent=0pt\hngind=\parindent \advance\hngind by 2.3em%
        \vskip0pt\hangafter=0\global\hangindent\hngind%
        {$\hbox to2.3em{#1\hskip.5em\hfill}$}}}%
\def\iitem{\global\advance\randomct by 1%
    \edef\numtoks{\the\randomct}%
    \idoitem{(\numtoks)}%
    \futurelet\testchar\maybeoptioniitem} 
\def\maybeoptioniitem{\ifx[\testchar\let\next\optioniitem%
	\else\let\next\nooptioniitem\fi\next}
\def\optioniitem[#1]{%
  \htmlanchor{#1}{\makeref{\noexpand#1}{Part}{\numtoks}}%
  \immediate\write\isauxout{\noexpand\forward{\noexpand#1}{}{(\numtoks)}}%
  \draftlabel{#1}}
\def\nooptioniitem{}

\let\ritem=\iitem

\def\citem{\global\advance\randomct by 1%
    \edef\numtoks{\the\randomct}%
    \idoitem{(\numtoks)}%
    \futurelet\testchar\maybeoptioncitem} 
\def\maybeoptioncitem{\ifx[\testchar\let\next\optioncitem%
	\else\let\next\nooptioniitem\fi\next}
\def\optioncitem[#1]{%
 \htmlanchor{#1}{\makeref{\noexpand#1}{Condition}{(\numtoks)}}%
 \immediate\write\isauxout{\noexpand\forward{\noexpand#1}{Condition}{(\numtoks)}}%
 \draftlabel{#1}%
 }

\def\idoiteitem#1{{\parindent=0pt\hngind=\parindent \advance\hngind by 3.9em%
        \vskip0pt\hangafter=0\global\hangindent\hngind%
        {$\hbox to3.9em{\hfill#1\ }$}}}
\def\idoitemitem#1{{\parindent=0pt\hngind=\parindent \advance\hngind by 4.9em%
        \vskip0pt\hangafter=0\global\hangindent\hngind%
        {$\hbox to4.9em{\hfill#1\ }$}}}
\def\idoitemitemm#1{{\parindent=0pt\hngind=\parindent \advance\hngind by 6.9em%
        \vskip0pt\hangafter=0\global\hangindent\hngind%
        {$\hbox to6.9em{\hfill#1\ }$}}}



\newdimen\thmskip \thmskip=0.7ex
	
\def\thmalready#1{\removelastskip\vskip\thmskip%
	\global\randomct=0%
	\noindent{\bf #1:}\ \ %
	\bgroup \it%
	\abovedisplayskip=4pt\belowdisplayskip=3pt%
	\parskip=0pt%
	}


\newdimen\exerciseindent \exerciseindent=4.5em
\newdimen\exerciseitemindent \exerciseitemindent=1.5em
\newdimen\exitindwk \exitindwk=3.5em
\def\exitem{\global\advance\randomct by 1%
	\par\indent\hangindent\exitindwk \hskip\exerciseitemindent%
	\llap{(\romannumeral\randomct)\enspace}\ignorespaces}
\def\stexitem{\global\advance\randomct by 1%
	\par\indent\hangindent\exitindwk \hskip\exerciseitemindent%
	\llap{(\romannumeral\randomct)\hskip-.5ex*\enspace}\ignorespaces}


\def\solution#1{\removelastskip\vskip 5pt%
	\global\randomct=0%
	\parindent=\exerciseindent%
	\par\hangindent\exerciseindent\indent\llap{\hbox to \exerciseindent{%
	  \bf \unskip\ref{#1}}}}

\def\solution#1{\removelastskip\vskip 5pt\par\noindent \ref{#1}}


\newdimen\wddim
\def\edinsert#1{\setbox0=\hbox{#1}\wddim=\wd0\divide\wddim by 2%
\smash{\hskip-\wddim\raise15pt\hbox to 0pt{\Red{$\underbrace{\box0}$}}\hskip\wddim}}


\def\hpad#1#2#3{\hbox{\kern #1\hbox{#3}\kern #2}}
\def\vpad#1#2#3{\setbox0=\hbox{#3}\vbox{\kern #1\box0\kern #2}}




\def\stack#1#2#3{\vbox{\offinterlineskip
  \setbox2=\hbox{#2}
  \setbox3=\hbox{#3}
  \dimen0=\ifdim\wd2>\wd3\wd2\else\wd3\fi
  \hbox to \dimen0{\hss\box2\hss}
  \kern #1
  \hbox to \dimen0{\hss\box3\hss}}}


\def\hexp#1{%
  \setbox0=\hbox{${}^{#1}$}%
  \hbox to .5\wd0{\box0\hss}}

\def\hsub#1{%
  \setbox0=\hbox{${}_{#1}$}%
  \hbox to .5\wd0{\box0\hss}}



\def\bmatrix#1#2{{\left(\vcenter{\halign
  {&\kern#1\hfil$##\mathstrut$\kern#1\cr#2}}\right)}}

\def\rightarrowmat#1#2#3{
  \setbox1=\hbox{\small\kern#2$\bmatrix{#1}{#3}$\kern#2}
  \,\vbox{\offinterlineskip\hbox to\wd1{\hfil\copy1\hfil}
    \kern 3pt\hbox to\wd1{\rightarrowfill}}\,}

\def\leftarrowmat#1#2#3{
  \setbox1=\hbox{\small\kern#2$\bmatrix{#1}{#3}$\kern#2}
  \,\vbox{\offinterlineskip\hbox to\wd1{\hfil\copy1\hfil}
    \kern 3pt\hbox to\wd1{\leftarrowfill}}\,}

\def\rightarrowbox#1#2{
  \setbox1=\hbox{\kern#1\hbox{\small #2}\kern#1}
  \,\vbox{\offinterlineskip\hbox to\wd1{\hfil\copy1\hfil}
    \kern 3pt\hbox to\wd1{\rightarrowfill}}\,}

\def\leftarrowbox#1#2{
  \setbox1=\hbox{\kern#1\hbox{\small #2}\kern#1}
  \,\vbox{\offinterlineskip\hbox to\wd1{\hfil\copy1\hfil}
    \kern 3pt\hbox to\wd1{\leftarrowfill}}\,}


\newcount\countA
\newcount\countB
\newcount\countC

\def\monthname{\begingroup
  \ifcase\number\month
    \or January\or February\or March\or April\or May\or June\or
    July\or August\or September\or October\or November\or December\fi
\endgroup}

\def\dayname{\begingroup
  \countA=\number\day
  \countB=\number\year
  \advance\countA by 0 
  \advance\countA by \ifcase\month\or
    0\or 31\or 59\or 90\or 120\or 151\or
    181\or 212\or 243\or 273\or 304\or 334\fi
  \advance\countB by -1995
  \multiply\countB by 365
  \advance\countA by \countB
  \countB=\countA
  \divide\countB by 7
  \multiply\countB by 7
  \advance\countA by -\countB
  \advance\countA by 1
  \ifcase\countA\or Sunday\or Monday\or Tuesday\or Wednesday\or
    Thursday\or Friday\or Saturday\fi
\endgroup}

\def\timename{\begingroup
   \countA = \time
   \divide\countA by 60
   \countB = \countA
   \countC = \time
   \multiply\countA by 60
   \advance\countC by -\countA
   \ifnum\countC<10\toks1={0}\else\toks1={}\fi
   \ifnum\countB<12 \toks0={\sevenrm AM}
     \else\toks0={\sevenrm PM}\advance\countB by -12\fi
   \relax\ifnum\countB=0\countB=12\fi
   \hbox{\the\countB:\the\toks1 \the\countC \thinspace \the\toks0}
\endgroup}

\def\timestamp{\dayname, \the\day\ \monthname\ \the\year, \timename}


\def\enma#1{{\ifmmode#1\else$#1$\fi}}

\def\mathbb#1{{\bbold #1}}
\def\mathbf#1{{\bf #1}}




%

%

%

%
%


\let\text\hbox
\def\overset#1\to#2{\ \mathop{\buildrel #1 \over #2}\ }
\def\underset#1\to#2{{\mathop{\buildrel #1 \over #2}}}

\def\myoperator{NOOP}
\def\mymathopsp{\futurelet\testchar\maybesubscriptop}
\def\maybesubscriptop{\ifx_\testchar \let\next\subscriptop
	\else \let\next\nosubscriptop \fi \next}
\def\subscriptop_#1{{\rm \myoperator}_{#1}}
\def\nosubscriptop{\futurelet\testchar\maybeparenop}
\def\maybeparenop{{\mathop{\rm \myoperator}\nolimits%
	\ifx(\testchar \relax\else\,\fi}}


\def\notdiv{\hbox{$\not\hbox to -0.2ex{$|$}$}\hskip0.5em}
\def\lt{{\rm lt}}
\def\lt{{\mathop{\rm lt}\nolimits}\,}
\def\lt{{\mathop{\rm lt}\nolimits}}

\def\sgn{\mathop{\rm sgn}\nolimits}


\def\lcm{\mathop{\rm lcm}\nolimits}


\def\lcm{\gdef\myoperator{lcm}\mymathopsp\!}


\def\uspace{\hbox{$\underline{\hbox to 1ex{\ \hfil}}$}}

\def\lineunderblank{\underline{\ \hbox to 0.5ex{$\vphantom{a}$\hfill\ }}}
\def\lineaboveblank{\overline{\ \hbox to 0.5ex{$\vphantom{a}$\hfill\ }}}
\def\blank#1{\hbox{$\underline{\hbox to #1em{\ \hfill}}$}}


\def\qedbox{\hbox{\vbox{\hrule\hbox{\vrule\kern3pt\vbox{\kern6pt}\kern3pt\vrule}\hrule}}}
\def\qed{\unskip\hfill\qedbox\vskip\thmskip} 
\def\eqed{\eqno\hbox{\quad\qedbox}} 



\def\xspace{\futurelet\lettoken\doxspace}
\def\doxspace{%
  \ifx\lettoken\bgroup\else%
  \ifx\lettoken\egroup\else%
  \ifx\lettoken\/\else%
  \ifx\lettoken\ \else%
  \ifx\lettoken~\else%
  \ifx\lettoken.\else%
  \ifx\lettoken!\else%
  \ifx\lettoken,\else%
  \ifx\lettoken:\else%
  \ifx\lettoken;\else%
  \ifx\lettoken?\else%
  \ifx\lettoken/\else%
  \ifx\lettoken'\else%
  \ifx\lettoken)\else%
  \ifx\lettoken-\else%
  \ifx\lettoken\space\else%
   \space%
   \fi\fi\fi\fi\fi\fi\fi\fi\fi\fi\fi\fi\fi\fi\fi\fi}


\def\today{\ifcase\month\or January\or February\or March\or
  April\or May\or June\or July\or August\or September\or
  October\or November\or December\fi
  \space\number\day, \number\year}\relax


\newwrite\isauxout
\openin1\jobname.aux
\ifeof1\message{No file \jobname.aux}
       \else\closein1\relax\input\jobname.aux
       \fi
\immediate\openout\isauxout=\jobname.aux
\immediate\write\isauxout{\relax}

\newwrite\iscontout
\immediate\openout\iscontout=\jobname.cont
\immediate\write\iscontout{\relax}
\def\contlist#1#2#3{\vskip1pt\noindent \ref{#1}: #2\hfill#3}
\def\contnosectlist#1#2#3{\vskip1pt\noindent #2\hfill#3}




\irenasize

\hyperstrue
\hypersfalse
\drafttrue
\draftfalse
\colortrue
\colorfalse
\ifcolor\input colordvi\else\input blackdvi\fi%
\ifdraft\footline={\hfill page \folio -- \today\hfill}\else
\footline={\hfill\folio\hfill}\fi

\def\switch{{\mathop{\rm s}\nolimits}}

\def\G{{G}}
\def\Gt{{\widetilde{G}}}

\def\I{{I}}
\def\It{{\widetilde{I}}}
\def\J{{J}}
\def\Jt{{\widetilde{J}}}
\def\Au_#1{{A^u_{#1}}}
\def\Gu_#1{{G^u_{#1}}}
\def\Fu_#1{{F^u_{#1}}}
\def\CI#1{{\I^{\langle #1 \rangle}}}
\def\CIt#1{{\It^{\langle #1 \rangle}}}
\def\CItS#1#2{{\It^{\langle #1 \rangle}_{#2}}}
\def\CJ#1{{\J^{\langle #1 \rangle}}}
\def\ChJ#1{{\widehat{\J}^{\langle #1 \rangle}}}
\def\CcJ#1{{\check{\J}^{\langle #1 \rangle}}}
\def\CJt#1{{\Jt^{\langle #1 \rangle}}}
\def\CJtS#1#2{{\Jt^{\langle #1 \rangle}_{#2}}}
\def\CG#1#2{{\Gt^{\langle #1 \rangle}_{#2}}}
\def\CVar#1#2{{\rm Var}^{\langle #1 \rangle}_{#2}}
\def\CQ#1#2{{\rm Q}^{\langle #1 \rangle}_{#2}}
\def\CQa#1#2{{\rm Q}^{\langle #1 \rangle}(#2)}
\def\CP#1#2{{\rm P}^{\langle #1 \rangle}_{#2}}
\def\pl{{\rm pl}}
\def\bpl{{\bf pl}}
\def\bipl{{\bf ipl}}
\def\ipl{{\rm ipl}}
\def\ipl{\overline {pl}}
\def\bipl{{\bf \overline {pl}}}
\def\bD{{D}}
\def\bbD{{D}}

\def\listmonoms{{\bf M}}

\def\pairslistvarsTswrad{{\bf T_3}} 
\def\pairslistvarsTswrad{{\bf T}}

\bigskip
\centerline{\titlefont $\bf 2 \times 2$ permanental ideals of hypermatrices}
\bigskip
\centerline{\titlefont Julia Porcino and Irena Swanson}
\bigskip

\bgroup
\baselineskip=11pt\narrower\narrower
\noindent
{\bf Abstract.}
We study the structure of ideals generated
by some classes of $2 \times 2$ permanents of hypermatrices.
This generalizes~\cite{LS} on $2 \times 2$ permanental ideal
of generic matrices.
We compare the obtained structure to that
of the corresponding determinantal ideals in \cite{ST}:
while the notion of $t$-switchability introduced in \cite{ST}
plays a role for both permanental and determinantal ideals,
the permanents require further restrictions,
which in general increases the number of minimal primes.
In the last two section we examine a few related classes of permanental ideals.

\egroup 

\section{intro}{Introduction}

Determinants are ubiquitous in mathematics;
if in the Laplace expansion of a determinant we replace all minus signs
with plus signs,
the result is the {\it permanent} of the matrix.
The starting motivation for this work is the paper~\cite{ST}
on classes of $2 \times 2$ determinantal ideals of certain hypermatrices
that arise in models of conditional independence;
here we study the analogous $2 \times 2$ permanental ideals.

The computation of determinants can be done in polynomial time
because of Gaussian elimination,
whereas the computation of permanents is \#-P-complete (\cite{V}).
In another comparison,
the ideal generated by $r \times r$ determinants of
a generic $m \times n$ matrix is a prime ideal with many nice properties (\cite{BV}),
but the ideal generated by $2 \times 2$ permanents has
${m \choose 2} {n \choose 2} + m + n$ minimal components
and one embedded component if $m, n \ge 3$ (\cite{LS}),
and the structure of ideals generated by $3 \times 3$ or larger permanents
is even more complicated and not completely understood (\cite{Ki}).
This paper puts the structure of the $2 \times 2$ permanental ideals,
as described in \cite{LS},
into more general and new perspective.
See \ref{exadmtwodiml} for how \cite{LS} can be viewed in this context.

Conditional independence ideals in algebraic statistics
are ideals generated by some $2 \times 2$ determinants
of generic hypermatrices.
Their structures have been studied recently by
Fink~\cite{Fink}, 
Herzog--Hibi--Hreinsdottir--Kahle--Rauh~\cite{HHHKR},
Ohtani~\cite{O},  
Swanson--Taylor~\cite{ST},
and Ah--Rauh~\cite{AR}.
In this paper we study the structure
of ideals $\CJ t$ of
analogous $2 \times 2$ permanents of generic hypermatrices
that depend on a parameter $t$
(and the size of the hypermatrix, of course).
The main result is \ref{thmCJminprimes}
which gives a combinatorial desription of all the prime ideals minimal over $\CJ t$.
Not surprisingly,
just as in~\cite{LS},
the number of minimal prime ideals
is in general larger over these permanental ideals
than over the corresponding determinantal ideals.
We give an explicit set of generators of the minimal primes
and give their Gr\"obner bases (\ref{thmGB}).
\ref{sectexamples} shows many concrete examples.

In \refs{sectjuliahat} and \refn{sectjuliacheck}
we present the structure of related permanental ideals,
those generated by certain ``diagonal" permanents,
and by certan diagonal and slice permanents.

Many of the methods in this paper are similar to those in~\cite{ST},
but for permanental ideals one has to keep track of the signs,
which adds not only notational difficulty but also changes many results.
In particular,
it increases the number of minimal primes significantly.

This is an extension of the first author's senior thesis at Reed College, 2011,
under the second author's supervision.

\section{setup}{Set-up}

We fix positive integers $n$ and $r_1, \ldots, r_n$.
Let $R$ be the polynomial ring in $r_1 \cdots r_n$ variables over
a field $k$.
The variables will be written with lower case $x$ with $n$-place subscripts,
with the $i$th place in the subscript ranging from $1$ through $r_i$.
We arrange these variables into a generic
$r_1 \times \cdots \times r_n$ hypermatrix
that will be fixed throughout.

Throughout the paper $[r]$ denotes the set $\{1, 2, \ldots, r\}$.
Thus, for example,
the ring we use is $R = k[x_a : a \in [r_1] \times \cdots \times [r_n]]$.
Throughout $N = [r_1] \times \cdots \times [r_n]$.

Let $L \subseteq [n]$.
For $a, b \in N$ define the {\bf switch} function $\switch(L,a,b)$
that switches the $L$-entries of $a$ into $b$:
$\switch(L,a,b)$ is an element of $N$ whose $i$th component is
$$
\switch(L,a,b)_i = \cases{
b_i, & if $i \in L$; \cr
a_i, & otherwise. \cr
}
$$
If $L = \{j\}$,
we simply write 
$\switch(L,a,b) = \switch(j,a,b)$.
For any $a$ and $b$ in $N$ we define
{\bf the distance $\bf d$} between them to be
${\bf d(a,b)} = \#\{i: a_i \not = b_i\}$. 
Note that $d(a,b) = d(\switch(L,a,b), \switch(L,b,a))$.

For any subsets $K \subseteq L \subseteq [n]$
and any $t \in [n]$
we define:
$$
\eqalignno{
f_{K,a,b} &= x_a x_b - x_{\switch(K,a,b)} x_{\switch(K,b,a)}, \cr
g_{K,a,b} &= x_a x_b + x_{\switch(K,a,b)} x_{\switch(K,b,a)}, \cr
\G_{L,K} &= \{g_{K,a,b} : a, b \in N, \{j: a_j \not = b_j\} = L\}. \cr
}
$$
The elements $g_{K,a,b}$ above are (generalized) permanents,
and the $f_{K,a,b}$ are (generalized) determinants.
Whenever $K = \{i\}$,
we write $i$ instead of $\{i\}$,
such as $f_{i,a,b}$, $g_{i,a,b}$, $G_{L,i}$.
When $d(a,b) = 2$ and $a_i \not = b_i$,
we call $g_{i,a,b}$ a {\bf slice permanent}
and $f_{i,a,b}$ a {\bf slice determinant}.
For any $t \in [n]$ we define:
$$
\eqalignno{
\CI{t} &= (f_{i,a,b} : a, b \in N, d(a,b) = 2,  i \in [t], a_i \not = b_i), \cr
\CJ{t}&= (g_{i,a,b} : a, b \in N, d(a,b) = 2,  i \in [t], a_i \not = b_i). \cr
}
$$
We examine the minimal primes over $\CJ{t}$,
whereas \cite{ST} did so for $\CI{t}$.
The structure of ideals
generated by classes of slice permanents
turns out to be different from the structure
of the corresponding ideals generated by slice determinants~\cite{ST}.
A special case of this difference was already demonstrated in~\cite{LS}
when $n = 2$.
In order to have $\CJ t$ different from $\CI t$,
we assume throughout
that the characteristic of the underlying field is not~$2$.

\section{sectlemmas}{Lemmas for induction arguments}

It is proved in~\cite{ST} that
for $a, a_1, b \in N$ and $i \in [n]$,
$x_b f_{i,a,a_1} - x_{a_1}f_{i,a,b}
=
x_{\switch(i,b,a)} f_{i,a_1,\switch(i,a,b)}
- x_{\switch(i,a,a_1)}f_{i,\switch(i,a_1,a),b}$.
We prove the almost analogous result for permanental ideals:

\lemma[lmInductIJ]
For $a, a_1, b \in N$ and $i \in [n]$,
$$
x_b f_{i,a,a_1} - x_{a_1}f_{i,a,b}
=
x_{\switch(i,b,a)} g_{i,a_1,\switch(i,a,b)}
- x_{\switch(i,a,a_1)}g_{i,\switch(i,a_1,a),b}.
$$
In particular,
if $a$ and $a_1$ differ only in the $j$th component,
and $j \not = i$,
and $b_i \not = (a_1)_i$,
then
$$
x_a g_{i,a_1,b}
- x_{a_1}f_{i,a,b}
=
x_{\switch(i,b,a)} g_{i,a_1,\switch(i,a,b)}
\in (G_{\{i,j\},\{i\}}).
$$
\endb

\proof
Just as in~\cite{ST},
the calculation of the first part is straightforward:
$$
\eqalignno{
x_b f_{i,a,a_1} -
&x_{a_1}f_{i,a,b} + x_{\switch(i,a,a_1)}g_{i,\switch(i,a_1,a),b} \cr
&=
x_b x_a x_{a_1}
- x_b x_{\switch(i,a,a_1)} x_{\switch(i,a_1,a)}
- x_{a_1} x_a x_b
+ x_{a_1} x_{\switch(i,a,b)} x_{\switch(i,b,a)} \cr
&\hskip4em
+ x_{\switch(i,a,a_1)} x_{\switch(i,a_1,a)} x_b
+ x_{\switch(i,a,a_1)}
x_{\switch(i,\switch(i,a_1,a),b)}
x_{\switch(i,b,\switch(i,a_1,a))}
\cr
&=
x_{a_1} x_{\switch(i,a,b)} x_{\switch(i,b,a)}
+ x_{\switch(i,a,a_1)}
x_{\switch(i,a_1,b)}
x_{\switch(i,b,a)}
\cr
&= x_{\switch(i,b,a)} (x_{a_1} x_{\switch(i,a,b)}
+
x_{\switch(i,a,a_1)} x_{\switch(i,a_1,b)})  \cr
&= x_{\switch(i,b,a)} g_{i,a_1,\switch(i,a,b)}.
}
$$
If $a$ and $a_1$ differ only in the $j$th position with $j \not = i$,
then $f_{i,a,a_1} = 0$ for all $i$,
$\switch(i,a_1,a) = a_1$,
$\switch(i,a,a_1) = a$,
and $a_1$ and $\switch(i,a,b)$ differ at most in the two components
$i$ and $j$,
and the rest follows.
\qed

The following lemma
is via induction an immediate corollary of~\ref{lmInductIJ}:

\lemma[lmvarsmult]
Let $i$ be a positive integer,
and let $a_0, a_1, \ldots, a_k, b \in N$.
Suppose that the $i$th component of $b$
differs from the $i$th components of $a_1, \ldots, a_k$,
and that for all $j = 1, \ldots, k$,
$a_{j-1}$ and $a_j$ differ exactly in component $l_j \not = i$.
Then modulo $\sum_{j=1}^k (\G_{\{i,l_j\},\{i\}})$,
$$
\eqalignno{
x_{a_1} x_{a_2} \cdots x_{a_k} f_{i,a_0,b}
&\equiv \cases{
x_{a_0} x_{a_1} \cdots x_{a_{k-1}} g_{i,a_k,b}, & if $k$ is odd; \cr
x_{a_0} x_{a_1} \cdots x_{a_{k-1}} f_{i,a_k,b}, & if $k$ is even; \cr
} \cr
x_{a_1} x_{a_2} \cdots x_{a_k} g_{i,a_0,b}
&\equiv \cases{
x_{a_0} x_{a_1} \cdots x_{a_{k-1}} f_{i,a_k,b}, & if $k$ is odd; \cr
x_{a_0} x_{a_1} \cdots x_{a_{k-1}} g_{i,a_k,b}, & if $k$ is even. \cr
} \cr
}
$$
In particular,
if $d(a_k,b) = 2$ and $a_k$ and $b$ differ in positions $i$ and $l_0 \not = i$,
then $\sum_{j=0}^k (\G_{\{i,l_j\},\{i\}})$ contains
$x_{a_1} x_{a_2} \cdots x_{a_k} f_{i,a_0,b}$ if $k$ is odd,
and it contains
$x_{a_1} x_{a_2} \cdots x_{a_k} g_{i,a_0,b}$ if $k$ is even.
Also,
if $d(a_k,b) = 1$,
then $\sum_{j=0}^k (\G_{\{i,l_j\},\{i\}})$ contains
$x_{a_1} x_{a_2} \cdots x_{a_k} f_{i,a_0,b}$ if $k$ is even,
and it contains
$x_{a_1} x_{a_2} \cdots x_{a_k} g_{i,a_0,b}$ if $k$ is odd.
\qed
\endb

\section{sectadmgraphla}{Switchable and signed sets}

Just as in~\cite{ST},
we need $t$-switchability and connectedness:

\defin[defswitchable]
A subset $S$ of $N$ is {\bf $\bf t$-switchable}
if for all $a, b \in S$ with $d(a,b) = 2$,
and all $i \in [t]$,
we have that $\switch(i,a,b), \switch(i,b,a) \in S$.

If a set $S$ is $t$-switchable,
define $a, b \in S$ to be {\bf connected}
if there exist $c_0 = a, c_1, \ldots, c_{k-1}, c_k = b$ in $S$
such that for all $i$,
$d(c_i,c_{i+1}) = 1$.
We call $c_0 = a, c_1, \ldots, c_{k-1}, c_k = b$ as above
{\bf a path from $\bf a$ to $\bf b$}.
We call $k$ {\bf the length of the path}.
The smallest length of a path from $a$ to $b$
is called the {\bf path length from $\bf a$ to $\bf b$}
and will be denoted as $\bf \bpl_S(a,b)$.
We define $\bf \bipl_S(a,b)$ to be $\pl_S(a,b) - 1$.
\endb

\remark[rmktadmswitcht]
Clearly connectedness is an equivalence relation.
By Lemma 3.5 in~\cite{ST},
for any $a, b$ that are connected in some $t$-admissible set
and for any subset $K \subseteq [t]$,
$\switch(K,a,b), \switch(K,b,a)$ are in $S$
and connected to $a$ and $b$ in $S$.

For a given $a, b$ that are connected in $S$,
it may happen that there are paths of different lengths between them.
Furthermore,
there can be paths of lengths of different parity between them,
such as if $S = [3] \times [2] \times [2]$:
$(1,1,1), (2,1,1), (3,1,1)$
and $(1,1,1), (3,1,1)$ are both paths from $(1,1,1)$ to $(3,1,1)$.

\prop[proppathparity]
Let $S$ be a $t$-switchable set,
and let $a, b \in S$ be connected with $d(a,b) \ge 2$ and $a_i \not = b_i$
for some $i \in [t]$.
Suppose that there are paths from $a$ to $b$ of different parity.
Then
$\sum_{j \not = i} (\G_{\{i,j\},\{i\}})$ contains monomials
all of whose subscripts are elements of~$S$.
\endb

\proof
Let $c_0 = a, c_1, \ldots, c_{k-1}, c_k = b$ be a path from $a$ to $b$.
Since $d(a,b) \ge 2$,
it follows that $k \ge 2$.

Let $j$ be least such that the $i$th entry of $c_j$ is not $a_i$.
Necessarily $j > 0$.
Suppose that $j \le k-2$.
By assumption,
$c_{j-1}$ and $c_j$ differ precisely in the $i$th component.
Suppose that for all $l = j-1, \ldots, k-1$,
$c_l$ and $c_{l+1}$ differ exactly in the $i$th component.
Then we can cut an even number of elements $c_j, \ldots, c_{k-1}$
from the given path to still get a path from $a$ to $b$
with the same parity but such that the designated $j$ is at least $k-1$.
Now assume otherwise:
then there exists an integer $l > j$
such that $c_{l-1}$ and $c_l$ differ in component $j \not = i$.
We choose $l$ to be the smallest such integer.
Then the part $c_{j-1}, c_j, \ldots, c_l$ of the given path
can by \ref{rmktadmswitcht} be replaced by
$c_{j-1}, \switch(i,c_l,c_{j-1}),
\switch(i,c_l,c_j),
\switch(i,c_l,c_{j+1}),
\ldots,
\switch(i,c_l,c_{l-1}) = c_l$,
which has the same length,
and the designated $j$ is increased by $1$.
By repeating these steps we may assume that the $i$th entries
in $c_0, c_1, \ldots, c_{k-2}$ all equal $a_i$.

In particular,
$1 \le d(c_{k-2}, b) \le 2$.


First suppose that $d(c_{k-2}, b) = 1$.
Then $c_{k-2}, c_{k-1}, b$ have distinct $i$th entries,
and only differ in the $i$th entries.
Since $d(a,b) \ge 2$,
necessarily $k > 2$,
and so $c_{k-3}$ differs from $c_{k-2}, c_{k-1}, b$
in the $j$th entry for some $j \not = i$,
and differs from $c_{k-1}, b$ also in the $i$th entry.
Then
$$
\left[
\matrix{
x_{c_{k-3}} & x_{\switch(i, c_{k-3},c_{k-1})} & x_{\switch(i, c_{k-3},b)} \cr
x_{c_{k-2}} & x_{\switch(i, c_{k-2},c_{k-1})} = x_{c_{k-1}} & x_{\switch(i, c_{k-2},b)} = x_b \cr
}
\right]
$$
is a submatrix of the hypermatrix,
and all of its $2 \times 2$ permanents are in $(\G_{\{i,j\},\{i\}})$.
Let $p_{u,v}$ be the permanent of columns $u,v$ in the matrix above.
Then $(\G_{\{i,j\},\{i\}})$ contains
$$
x_{c_{k-3}} p_{2,3} - x_{\switch(i, c_{k-3},c_{k-1})} p_{1,3}
+ x_{\switch(i, c_{k-3},b)} p_{1,2}
= 2
x_{c_{k-3}} x_{c_{k-1}} x_b,
$$
and as the characteristic is not $2$,
$(\G_{\{i,j\},\{i\}})$ contains a monomial
all of whose subscripts are in $S$.

So we may assume that $d(c_{k-2}, b) = 2$
for all paths from $a$ to $b$.
Similarly,
we may assume that if $d_0 = a, d_1, \ldots, d_{l-1}, d_l = b$
is a path from $a$ to $b$
with $l$ of different parity from $k$,
then the $i$th components of $d_0, d_1, \ldots, d_{l-2}$ are all $a_i$
and $d(d_{l-2},b) = 2$.
Say $k$ is even and $l$ is odd.
Then by \ref{lmvarsmult},
$x_{c_1} \cdots x_{c_{k-2}} g_{i,a,b},
\quad
x_{d_1} \cdots x_{d_{l-2}} f_{i,a,b} \in
\sum_{j\not=i} (\G_{\{i,j\},\{i\}})$.
In particular,
$$
\lcm(x_{c_1} \cdots x_{c_{k-2}}, x_{d_1} \cdots x_{d_{l-2}})
(g_{i,a,b}, f_{i,a,b}) \subseteq \sum_{j\not=i} (\G_{\{i,j\},\{i\}}).
$$
But $(g_{i,a,b}, f_{i,a,b}) 
= (x_a x_b, x_{\switch(i,a,b)} x_{\switch(i,b,a)})$
as ideals,
so that $\sum (\G_{\{i,j\},\{i\}})$ contains the monomials
$$
\lcm(x_{c_1} \cdots x_{c_{k-2}}, x_{d_1} \cdots x_{d_{l-2}}) x_a x_b
\hbox{ and }
\lcm(x_{c_1} \cdots x_{c_{k-2}}, x_{d_1} \cdots x_{d_{l-2}})
x_{\switch(i,a,b)} x_{\switch(i,b,a)}.  \eqed
$$

For our analysis of prime ideals that are minimal over
the permanental ideal $\CJ{t}$,
we deal with the parity question
with the following further restriction on switchable sets:

\defin[deftadm]
A $t$-switchable set $S$ is called {\bf $\bf t$-signed}
if for every equivalence class~$S_0$ with respect to the connected property,
one of the following conditions is satisfied:
\ritem
All elements of $S_0$ have the same first $t$-coordinates;
\ritem
Any two elements of $S_0$ differ at most in one component;
\ritem
The parity of path length between any two elements in $S_0$
is independent of the path.
\endb

\remark[rmkparityswitch]
Let $S$ be a $t$-signed set
and let $a, b, c \in S$ be connected
and in an equivalence class that satisfies property (3) in \ref{deftadm}.
Then for any $K \subseteq [t]$,
$\pl_S(a,b)$ and $\pl_S(\switch(K,a,b),\switch(K,b,a))$
have the same parity
because we can make the path from $\switch(K,a,b)$ to $\switch(K,b,a)$
pass first from $\switch(K,a,b)$ to $a$
in $\#\{i \in K: a_i \not = b_i\}$ steps,
then to $b$ in $\pl_S(a,b)$ steps,
and then to $\switch(K,b,a)$ in $\#\{i \in K: a_i \not = b_i\}$ steps.
Similarly,
$\pl_S(a,b) + \pl_S(b,c)$ and $\pl_S(a,c)$ have the same parity.

\lemma[lmtadm2]
Let $S$ be a $t$-signed set.
Let $S_0$ be an equivalence class in $S$ with respect to connectedness.
Then either for all $i \in [t]$,
$\{a_i : a \in S_0\}$ contains at most two elements
or else for all $a, b \in S_0$,
$d(a,b) \le 1$.
\endb

\proof
Suppose for contradiction that there are $a, b, a', b', c' \in S_0$
such that for some $i \in [t]$,
the $i$th components of $a', b', c'$ are all distinct
and such that $d(a,b) \ge 2$.
Thus $S_0$ fails conditions (1) and (2) in \ref{deftadm},
so $S_0$ must satisfy (3).
But then
$\switch(i,a,a'), \switch(i,a,b'), \switch(i,a,c')$
and $\switch(i,a,a'), \switch(i,a,c')$
are both paths from $\switch(i,a,a')$ to $\switch(i,a,c')$,
and they have different parities,
which contradicts the condition in (3).
%
\qed

One might be tempted to think that
for a $t$-signed equivalence class $S_0$,
if there exists $i \in [t]$ such that $\{a_i : a \in S_0\}$ has two elements,
then for all $j \in [n]$,
$\{a_j : a \in S_0\}$ has at most $2$ elements.
This need not be the case,
as demonstrated in \ref{ex223-1adm}~(4).

The following encapsulates and generalizes the $f_{K,a,b}$ and $g_{K,a,b}$,
and will be used to describe generators of prime ideals
minimal over $\CJ t$ (details in \ref{rmkhwelldef}):

\defin[defadmQ]
Let $S$ be $t$-signed.
For any connected $a, b \in S$,
and any $K \subseteq \{i \in [t]: a_i \not = b_i\}$,
define
$$
h_{S,K,a,b}
= x_a x_b
- (-1)^{\#K \cdot \ipl_S(a,b)} x_{\switch(K,a,b)}, x_{\switch(K,b,a)}.
$$
\endb

\remark[rmkhwelldef]
Suppose that for some $L \subseteq \{i \in [t]: a_i \not = b_i\}$,
$x_{\switch(K,a,b)} x_{\switch(K,b,a)}
= x_{\switch(L,a,b)} x_{\switch(L,b,a)}$.
We prove that then $h_{S,K,a,b} = h_{S,L,a,b}$.
If $\ipl_S(a,b)$ is even,
this is certainly true.
So we may assume that $\pl_S(a,b) = \ipl_S(a,b) + 1$ is even.
Without loss of generality $K \not = L$.
The assumption
$x_{\switch(K,a,b)}, x_{\switch(K,b,a)}
= x_{\switch(L,a,b)}, x_{\switch(L,b,a)}$
is only possible if $(a_{t+1}, \ldots, a_n) = (b_{t+1}, \ldots, b_n)$
or if $t = n$.
In either case,
$L = \{i : a_i \not = b_i\} \setminus K$,
and
$\{i : a_i \not = b_i\} \subseteq [t]$,
so that by \ref{rmktadmswitcht},
$\pl_S(a,b) = d(a,b) = \{i : a_i \not = b_i\}$.
But $\pl_S(a,b)$ is even,
so that $\#L$ and $\#K$ must have the same parity.
Thus $h_{S,K,a,b} = h_{S,L,a,b}$.

The following definitions will be applied mostly to $t$-signed $S$:

\defin[defadmPQ]
For any subset $S$ of $N$,
define
$$
\eqalignno{
\CVar{t}{S} &= (x_a : a \not \in S), \cr
\CJtS{t}{S} &= (h_{S,i,a,b} :
a, b \hbox{ are connected in $S$}, i \in [t], a_i \not = b_i), \cr
\CG{t}{S} &= \{h_{S,K,a,b} :
a, b \hbox{ are connected in $S$},
K \subseteq \{i \in [t]: a_i \not = b_i\}\}, \cr
\CQ{t}{S} &= \CVar{t}{S} + \CJtS{t}{S}. \cr
}
$$
\endb

Note that $t = n$ and $t = n-1$ give the same sets of ideals.
Thus in the sequel, and especially in \ref{sectGB},
we mostly talk about $t < n$.
For $t$-switchable $S$,
the analogous determinantal ideals $\CItS{t}{S}
= (f_{i,a,b} : i \in [t], a, b \hbox{ are connected in $S$})$
and
$\CP{t}{S} = \CVar{t}{S} + \CItS{t}{S}$
were proved to be prime ideals in \cite{ST}.
Much of what we prove in this paper mimics the proofs of~\cite{ST},
but with the added sign difficulty expressed through $t$-signedness.

\lemma[lmGtcont]
$\CG{t}{S} \subseteq \CJtS{t}{S}$.
\endb

\proof
Let $a, b \in S$
and let $K \subseteq \{i \in [t]: a_i \not = b_i\}$.
Write $K = \{k_1, \ldots, k_l\}$.
By \ref{rmktadmswitcht},
$\switch(\{k_1, \ldots, k_{i-1}\},a,b)$ and
$\switch(\{k_1, \ldots, k_{i-1}\},b,a))$ are connected,
and by \ref{rmkparityswitch},
$\pl_S(a,b) =
\pl_S(\switch(\{k_1, \ldots, k_{i-1}\},a,b),
\switch(\{k_1, \ldots, k_{i-1}\},b,a))$,
so that
$$
h_{S,K,a,b}
=
\sum_{i=1}^l
(-1)^{(i-1) \cdot \ipl_S(a,b)}
h_{S,k_i,\switch(\{k_1, \ldots, k_{i-1}\},a,b),
\switch(\{k_1, \ldots, k_{i-1}\},b,a)}.
\eqed
$$

\prop[propCJtinCQ]
Let $S$ be a $t$-switchable set (not necessarily $t$-signed).
Then $\CJ{t}~\subseteq~\CQ{t}{S}$.
\endb

\proof
Let $a, b \in N$ with $d(a,b) = 2$ and $a_i\neq b_i$ for some $i\in[t]$.
We need to prove that $h_{S,i,a,b} \in \CQ{t}{S}$.
If $a, b \in S$,
then by $t$-switchability of $S$,
$a$ and $b$ are connected,
so that $h_{S,i,a,b} \in \CJtS{t}{S} \subseteq \CQ{t}{S}$.
Similarly,
if $\switch(i,a,b), \switch(i,b,a) \in S$,
then
$h_{S,i,a,b} = \pm h_{S,i,\switch(i,a,b), \switch(i,b,a)} \in \CQ{t}{S}$,
So we may assume that $a \not \in S$
and that either
$\switch(i,a,b)$ or $\switch(i,b,a)$ is not in $S$.
But then $h_{S,i,a,b} \in \CVar{t}{S} \subseteq \CQ{t}{S}$.
\qed

\defin[defmaxtadm]
Let $S$ be $t$-signed.
We say that $S$ is {\bf maximal $\bf t$-signed}
if for all $t$-signed subsets $T$ of $N$ containing $S$,
$\CQ{t}{S} \subseteq \CQ{t}{T}$.
(So $\CVar{t}{S}$ contains $\CVar{t}{T}$,
that is why $S$ is called maximal.)
\endb

\section{sectexamples}{Some examples of switchable and signed sets}

To bring the abstract notion of $t$-switchable sets down to earth,
we now examine some concrete examples.
In \ref{sectminprimescont}
we characterize more generally
all minimal prime ideals over $\CJ t$
(they are $\CQ{t}{S}$ for maximal $t$-signed $S$),
but in examples below we simply state what the associated prime ideals are.
We used Macaulay2~\cite{GS} together with Kahle's binomial package~\cite{Ka}.

When the set $S$ has a long name,
we write $\CQa{t}{\!S}$ instead of $\CQ{t}{S}$.

First we provide perspective in light of $t$-signed sets
on the primary decomposition of $2 \times 2$ permanental ideals
of generic matrices as first determined without $t$-signed notion in~\cite{LS}.
The main result of~\cite{LS}
was to give the description of the minimal primes over
the ideal generated by $2 \times 2$ permanents of a generic $r_1 \times r_2$ matrix
in terms of three types primes;
but by the work in this paper,
the three types of primes are all simply the prime ideals
corresponding to $1$-signed sets:

\example[exadmtwodiml]
Let $N = [r_1] \times [r_2]$,
with $r_1, r_2 \ge 2$,
which represents an ordinary matrix.
For $t = 1, 2$,
the $t$-signed subsets of $N$ are identical,
and they are:
all $2 \times 2$,
$1 \times r$ and $r \times 1$ submatrices of~$N$.
If $r_1 = r_2 = 2$,
then the only maximal $t$-signed set is $N$.
If $r_1 = 2 < r_2$,
then the only maximal $t$-signed sets are
$2 \times 2$ and $1 \times r_2$ submatrices of~$N$,
and if $r_1, r_2 > 2$,
then the $t$-signed subsets of $N$ are $2 \times 2$,
$1 \times r_2$ and $r_1 \times 1$ submatrices of~$N$.
By \ref{thmCJminprimes},
these maximal $t$-signed sets correspond precisely
to the prime ideals minimal over $\CJ t$,
and this was first established in~\cite{LS}
without the vocabulary of signed sets.
There it was also established that when $r_1, r_2 > 2$,
there is exactly one embedded prime,
corresponding to the (non-maximal) $t$-signed set $\emptyset$.

\example[ex222notold1adm]
Let $N = [2] \times [2] \times [2]$.
Here we will think of $N$ as a cube.
The $1$-signed subsets of $N$,
together with their corresponding ideals, are as follows:
\ritem
$N$ (the whole cube);
$\CQ{1}{N} = \CJtS{t}{N}$
is generated by four slice permanental
and two diagonal determinantal generators.
\ritem
$S$ is a face of the cube with nonconstant first component;
$\CQa{1}{\hbox{face}}$ is generated by the permanent of this face
and the variables on the opposite face.
\ritem
$S$ cannot be a face of the cube perpendicular to the $x$-axis
as that is not even $1$-switchable;
\ritem
$S$ is an edge of the cube;
$\CQa{1}{\hbox{edge}}$ is generated by the variables not on this edge.
\ritem
two parallel edges of the cube with non-constant first component
that are not on the same face,
this is an example of a $1$-switchable set with two equivalence
classes with respect to connectedness;
$\CQa{1}{\hbox{two edges}}$ is generated by the variables not on these edges.
\ritem
a point;
$\CQa{1}{\hbox{point}}$ is generated by the other $7$ variables.
\ritem
two points of distance three
or two points of distance two with the same first component;
$\CQa{1}{\hbox{two points}}$ is generated by the variables not on these points.

Note that $\CQ{1}{N}$ is contained in $\CQa{1}{\hbox{face}}$,
$\CQa{1}{\hbox{point}}$, and $\CQa{1}{\hbox{two points}}$,
and that $\CQa{1}{\hbox{two edges}}$ is contained in $\CQa{1}{\hbox{edge}}$.
From this we read that the maximal $1$-signed sets
are: $N$, and two edges parallel to the $x$-axis that are not on the same face.
Note that there are two options for $S$ being formed by the parallel edges.
There are thus three maximal signed sets.
Macaulay2 \cite{GS} says that these $\CQa{t}{S}$
are all the minimal primes
(and even all the components).

\example[ex222-2adm]
Let $N = [2] \times [2] \times [2]$ as in the previous example.
The $2$-signed subsets of $N$ (which are the same as $3$-signed sets),
together with their corresponding ideals, are as follows:
\ritem
$N$ (the whole cube);
$\CQ{2}{N} = \CJtS{2}{N}$
is (redundantly) generated by one slice permanent for each of the six faces,
and three diagonal determinants for each of the opposite vertices.
\ritem
$S$ is a face of the cube;
$\CQa{2}{\hbox{face}}$ is generated by the permanent of this face
and the variables on the opposite face.
\ritem
an edge of the cube;
$\CQa{2}{\hbox{edge}}$ is generated by the variables not on this edge.
\ritem
a point;
$\CQa{2}{\hbox{point}}$ is generated by the remaining $7$ variables.
\ritem
two points of distance three;
$\CQa{2}{\hbox{two points}}$ is generated by the variables not on these points.

Note that the union of two parallel edges that are not on the same face
is not $2$-signed.
Here the maximal $2$-signed sets are:
$N$, and the two-opposite points.
There are thus five maximal signed sets,
and, indeed, Macaulay2 says that these are all the minimal components
(and even all the associated primes).

\example
With $N = [3] \times [2] \times [2]$,
the maximal $1$-signed sets are:
\ritem
For each $k \in [3]$,
$S$ is the set of variables whose first coordinate is not $k$,
so it looks like a $2 \times 2 \times 2$-block;
$\CQ{1}{S}$ is generated by the variables whose first coordinate is $k$,
and by the four slice permanents and two diagonal determinants on the block.
There are three such $S$.
\ritem
$S$ consists of the variables on two lines parallel to the $x$-axis,
the lines not lying on the same face of the cube;
$\CQ{1}{S}$ is generated by the variables not on these lines.
There are two such $S$.

Indeed,
Macaulay2 gives $3 + 2$ minimal primes,
but it also gives one embedded prime ideal consisting of all variables.
While the minimal components are all prime ideals themselves,
the embedded component in this case may be taken to be
$\CJ{1} + (x_a^3: a \in N)$.

\example[ex322-2adm]
With $N = [3] \times [2] \times [2]$ as in the previous example,
the maximal $2$-signed sets are:
\ritem
For each $i \in [3]$,
$S$ is the set of variables whose first coordinate is not $i$;
$\CQ{2}{S}$ is generated by the variables whose first coordinate is $i$,
and by the slice permanents and diagonal determinants
on the $2 \times 2 \times 2$-block.
There are three such $S$.
\ritem
$S$ is a line parallel to the $x$-axis;
$\CQ{2}{S}$ is generated by the $9$ variables not in $S$.
There are four such $S$.
\ritem
$S$ is the disjoint union of a point
and of two points on a line parallel to the $x$-axis,
with all distances between the one point and the other two points being $3$;
$\CQ{2}{S}$ is generated by the $9$ variables not in $S$.
There are $4 \cdot 3$ such $S$.

Indeed,
Macaulay2 gives $3 + 4 + 12 = 19$ minimal primes,
and it also gives one embedded prime ideal consisting of all variables.
The minimal components are all prime ideals,
and the embedded component in this case as well may be taken to be
$\CJ{2} + (x_a^3: a \in N)$.

\example[ex223-1adm]
With $N = [2] \times [2] \times [3]$,
by symmetry the maximal $2$-signed sets and primary decomposition
are as in \ref{ex322-2adm}.
(In contrast,
we do not have symmetry for $1$-signed sets.)
The maximal $1$-signed sets are as follows:
\ritem
$S$ consists of all variables that do not have a certain third coordinate
(so it is a $2 \times 2 \times 2$ subhypermatrix);
$\CQ{1}{S}$ is generated by the variables not in $S$
and by the diagonal determinants and the slice permanents of the block.
There are $3$ such $S$,
one for each of the three third coordinates.
\ritem
$S$ consists of variables on a face perpendicular to the $x$-axis;
$\CQ{1}{S}$ is generated by the variables not in $S$.
There are two such $S$.
\ritem
$S$ is a disjoint union of a line segment and of a $2 \times 2$ submatrix
on distinct faces perpendicular to the $y$-axis;
$\CQ{1}{S}$ is generated by one permanent of the $2 \times 2$-matrix
and by all variables not in $S$.
There are $2 \cdot 3$ such $S$.
\ritem
$S$ is the complement of the union of two lines parallel to the $x$-axis,
the two lines having distinct second and third coordinates.
This $S$ is a non-disjoint union of three $2 \times 2$ submatrices,
with one of the submatrices sharing edges with the other two.
With this visualization,
$\CQ{1}{S}$ is generated by all variables not in $S$,
by one permanent for each of the submatrices,
by one diagonal determinant on the elements of $S$
for each of the two pairs of adjacent $2 \times 2$ submatrices,
and by one extra diagonal permanent of the two disjoint edges of the
two $2 \times 2$ submatrices,
There are $3 \cdot 2$ such $S$.
(Note that this last $S$ is one equivalence class under connectedness;
there are two possible first coordinates for each element,
but the third coordinates can take on three different values
for various first two coordinates.)

Indeed,
Macaulay2 gives $17 = 3 + 2 + 6 + 6$ minimal primes,
and no embedded primes.
All the minimal components are prime.

\section{sectGB}{Primes, Gr\"obner bases and similar reductions}

The results here mimic those of \cite{ST}.
The generators of $\CJt t$ are (monic) binomials
that are, up to sign of the second coefficient in the binomial,
the same as the generators of $\CIt t$.
By Lemma~6.2 in~\cite{ST},
every monomial reduces modulo $\CIt t$ to a well-understood monomial;
thus by the nature of the generators of $\CJt t$,
every monomial reduces modulo $\CJt t$ to a well-understood monomial
times a less well-understood sign.
The bulk of this section is developing the machinery to keep track of this sign:
\ref{lmGB} introduces an absolute measure
for keeping track of what degree-three binomials are in $\CJtS{t}{S}$,
and \ref{lmGBcont}
proves that this measure behaves well
under any partial switches in one of the two monomials in the binomial.
The ultimate goal of this section is determining a Gr\"obner basis of $\CJt t$
and of proving that $\CJt t$ is a prime ideal.

Throughout this section we use the lexicographic order on the variables,
with variables sorted in the lexicographic order on their indices.
\def\blt{\hbox{\bf lt}}
By $\blt$ we denote the leading term.
Also,
throughout $S$ is a $t$-signed set,
and for slight brevity
we write $\pl_S(a,b), h_{S,K,a,b}$ without~``$S$".


For the purposes of discussion below we introduce the following notion:
we enumerate in the lexicographic order all possible elements of
$[r_{t+1}] \times \cdots \times [r_n]$,
with largest elements in the lexicographic order
getting the largest numeral.
Thus we can think of each element $a \in N$
as a $(t+1)$-tuple
in $[r_1] \times \cdots \times [r_t] \times [r_{t+1} \cdots r_n]$,
and in this notation without loss of generality $t = n-1 < n$.

\defin[deftilde]
When we think of $a$ as a $(t+1)$-tuple,
we denote $a$ as $\tilde a$.
\endb

\lemma[lmGB]
Let $S$ be a $t$-signed set,
and let $a,b,c,A,B,C \in S$ be mutually connected
such that $\widetilde a_{t+1} = \widetilde A_{t+1}$,
$\widetilde b_{t+1} = \widetilde B_{t+1}$,
$\widetilde c_{t+1} = \widetilde C_{t+1}$,
and for each $i = 1, \ldots, t$,
the list $a_i,b_i,c_i$ is up to order the same as the list $A_i,B_i,C_i$
(with same multiplicities).
Define
$$
\eqalignno{
K_1 &= \{i \in [t]: A_i = b_i \not = a_i\}, \cr
K_2 &= \{i \in [t]: A_i = c_i \not = a_i\} \setminus K_1, \cr
K_3 &= \{i \in K_1: B_i \not = a_i\} \cup
(\{i \in [t]: B_i \not = b_i\} \setminus K_1), \cr
p &= \#K_1 \cdot \ipl(a,b)
+ \#K_2 \cdot \ipl(\switch(K_1,a,b),c)
+ \#K_3 \cdot \ipl(\switch(K_1,b,a),\switch(K_2,c,a)). \cr
}
$$
Then $x_a x_b x_c - (-1)^p x_A x_B x_C \in \CJtS{t}{S}$.
\endb

\proof
$K_1$ and $K_2$ are disjoint by construction,
so that $\switch(K_2,c,a) = \switch(K_2,c,\switch(K_1,a,b))$.

Note that $\switch(K_2,\switch(K_1,a,b),c) = A$,
that $\switch(K_3,\switch(K_1,b,a),c) = B$,
and that $\switch(K_3,\switch(K_2,c,\switch(K_1,a,b)), \switch(K_1,b,a)) = C$.
Thus by \ref{lmGtcont},
the following elements are in $\CJtS{t}{S}$:
$$
\eqalignno{
&x_a x_b
- (-1)^{\#K_1 \cdot \ipl(a,b)} x_{\switch(K_1,a,b)} x_{\switch(K_1,b,a)}, \cr
&x_{\switch(K_1,a,b)} x_c
- (-1)^{\#K_2 \cdot \ipl(\switch(K_1,a,b),c)} x_A x_{\switch(K_2,c,\switch(K_1,a,b))}, \cr
&x_{\switch(K_1,b,a)} x_{\switch(K_2,c,\switch(K_1,a,b))}
- (-1)^{\#K_3 \cdot \ipl(\switch(K_1,b,a),\switch(K_2,c,\switch(K_1,a,b)))}
x_B x_C. \cr
}
$$
It follows that
$$
\eqalignno{
x_a &x_b x_c -
(-1)^{\#K_1 \cdot \ipl(a,b) + \#K_2 \cdot \ipl(\switch(K_1,a,b),c)
+ \#K_3 \cdot \ipl(\switch(K_1,b,a),\switch(K_2,c,\switch(K_1,a,b)))}
x_A x_B x_C \cr
&= \left(x_a x_b
- (-1)^{\#K_1 \cdot \ipl(a,b)} x_{\switch(K_1,a,b)} x_{\switch(K_1,b,a)}\right) x_c
\cr
&\hskip2em+
(-1)^{\#K_1 \cdot \ipl(a,b)} x_{\switch(K_1,b,a)}
\left(x_{\switch(K_1,a,b)} x_c
- (-1)^{\#K_2 \cdot \ipl(\switch(K_1,a,b),c)} x_A x_{\switch(K_2,c,\switch(K_1,a,b))}\right) \cr
&\hskip2em+
(-1)^{\#K_1 \cdot \ipl(a,b)
+ \#K_2 \cdot \ipl(\switch(K_1,b,a),c)} \cr
&\hskip4em \cdot
x_A \left(x_{\switch(K_1,b,a)} x_{\switch(K_2,c,\switch(K_1,a,b))}
- (-1)^{\#K_3 \cdot \ipl(\switch(K_1,b,a),\switch(K_2,c,\switch(K_1,a,b)))}
x_B x_C\right) \cr
}
$$
is in $\CJtS{t}{S}$.
\qed

Note that the order of the pairs $(a,A)$, $(b,B)$, $(c,C)$
in the lemma above affects the number $p$.
For this reason,
the proof of the next result requires many subcases,
but it leads to the eventual fact that the parity of $p$
is independent of the order of the pairs.


\lemma[lmGBcont]
With $a, b, c, A, B, C$ as in $\ref{lmGB}$,
we fix $A, B, C$,
and we allow $a, b, c$ to vary.
Thus $p$ is now a function of $a, b, c$.
Let $i \in [t]$,
and assume one of the following:
\ritem
$a_i \not = b_i$,
and
$a' = \switch(i,a,b)$, $b' = \switch(i,b,a)$, $c' = c$,
$r = \ipl(a,b)$.
\ritem
$a_i \not = c_i$,
and
$a' = \switch(i,a,c)$, $b' = b$, $c' = \switch(i,c,a)$,
$r = \ipl(a,c)$.
\ritem
$b_i \not = c_i$,
and
$a' = a$, $b' = \switch(i,b,c)$, $c' = \switch(i,c,b)$,
$r = \ipl(b,c)$.

Then $p(a',b',c') - p(a,b,c) + r$ is even.
\endb

\proof
By assumption and by \ref{rmktadmswitcht},
$a, b, c, a', b', c', A, B, C$ are all in the same
connected equivalence class.

%
To verify the evenness,
we need to recall the construction of \ref{lmGB}:
to get from $a,b,c$ to $A, B, C$,
we first switch $\#K_1$ entries from $b$ into $a$
and then $\#K_2$ entries from $c$ into $a$
to then have the modified $a$ equal to $A$;
the procedure is to choose $K_1$ as large as possible,
and then $K_2$ follows uniquely.
In the third step using $K_3$,
we switch the necessary entries of the new $b$ and the new $c$
to get the modified $b$ equal $B$,
and so necessarily the modified $c$ equals~$C$.
Let $K'_1, K'_2, K'_3$ be the corresponding sets
when we start with $a', b', c'$ in place of $a, b, c$.

Suppose that $a,b,c$ all differ at most in the $i$th component.
Then all $\ipl$ in the definitions of $p(a,b,c)$ and $p(a',b',c')$ are $0$,
and $r = 0$,
so that the lemma holds in this case.
So we may assume that at least two of $a,b,c$ 
differ not just in the $i$th component but also in some other $j$th component.
By \ref{lmtadm2} $a,b,c$,
the set $\{a_i,b_i,c_i\}$ contains exactly two elements.

Here are a few simplified expressions of $p(a,b,c)$ modulo $2$:
$$
\eqalignno{
&\ipl(a,b) = \pl(a,b) - 1, \cr
&\ipl(\switch(K_1,a,b),c)
=
\pl(\switch(K_1,a,b),c) - 1
\equiv
\pl(\switch(K_1,a,b),a) + \pl(a,c) - 1
\cr
&\hskip4em=
\#K_1 + \pl(a,c) - 1,
\cr
&\ipl(\switch(K_1,b,a),\switch(K_2,c,a))
=
\pl(\switch(K_1,b,a),\switch(K_2,c,a)) - 1
\cr
&\hskip4em\equiv
\pl(\switch(K_1,b,a),b) + \pl(b,c) + \pl(c,\switch(K_2,c,a)) - 1 \cr
&\hskip4em\equiv
\#K_1 + \pl(b,c) + \#K_2 - 1, \cr
&p(a,b,c) =
\#K_1 \cdot (\pl(a,b) - 1)
+ \#K_2 \cdot (\#K_1 + \pl(a,c) - 1) \cr
&\hskip8em
+ \#K_3 \cdot (\#K_1 + \pl(b,c) + \#K_2 - 1), \cr
}
$$
and there is a similar expressions for $p(a', b', c')$ in terms of
$K_1', K_2', K_3'$.

We analyze the cases separately.

(1)
As stated, by \ref{lmtadm2} we have two subcases:
$A_i = a_i \not = b_i$,
and $a_i \not = b_i = A_i$.
In the first subcase,
$i \in K_1' \setminus K_1$,
the only change in the construction of $A, B, C$
is that $K'_1 = K_1 \cup \{i\}$
(but $\ipl(a',b') = \ipl(a,b)$),
whence by \ref{lmGB} and the simplifications above,
modulo $2$,
$$
\eqalignno{
p(a',b',c') &- p(a,b,c) + r =
\#K_1' \cdot (\pl(a',b') - 1)
+ \#K_2' \cdot (\#K_1' + \pl(a',c') - 1) \cr
&\hskip4em
+ \#K_3' \cdot (\#K_1' + \pl(b',c') + \#K_2' - 1) \cr
&\hskip4em
- \#K_1 \cdot (\pl(a,b) - 1)
- \#K_2 \cdot (\#K_1 + \pl(a,c) - 1) \cr
&\hskip4em
- \#K_3 \cdot (\#K_1 + \pl(b,c) + \#K_2 - 1)
+ \ipl(a,b)
\cr
&\equiv (\#K_1 + 1) \cdot (\pl(a,b) - 1)
+ \#K_2 \cdot (\#K_1 + 1 + \pl(a',c) - 1) \cr
&\hskip4em
+ \#K_3 \cdot (\#K_1  + 1 + \pl(b',c) + \#K_2 - 1) \cr
&\hskip4em
- \#K_1 \cdot (\pl(a,b) - 1)
- \#K_2 \cdot (\#K_1 + \pl(a,c) - 1) \cr
&\hskip4em
- \#K_3 \cdot (\#K_1 + \pl(b,c) + \#K_2 - 1)
+ \pl(a,b) - 1
\cr
&\equiv (\#K_1 + 1) \cdot (\pl(a,b) - 1)
+ \#K_2 \cdot (\#K_1 + \pl(a',a) + \pl(a,c)) \cr
&\hskip4em
+ \#K_3 \cdot (\#K_1  + \pl(b',b) + \pl(b,c) + \#K_2) \cr
&\hskip4em
- \#K_1 \cdot (\pl(a,b) - 1)
- \#K_2 \cdot (\#K_1 + \pl(a,c) - 1) \cr
&\hskip4em
- \#K_3 \cdot (\#K_1 + \pl(b,c) + \#K_2 - 1)
+ \pl(a,b) - 1
\cr
&\equiv (\#K_1 + 1) \cdot (\pl(a,b) - 1)
- \#K_1 \cdot (\pl(a,b) - 1)
+ \pl(a,b) \cr
&\equiv \pl(a,b) - 1
+ \pl(a,b) - 1, \cr
}
$$
which is even.
In the second subcase ($a_i \not = b_i = A_i$),
$K_1' = K_1 \setminus \{i\}$
and there are no other changes,
so that the parity of
$p(a',b',c') - p(a,b,c) + r$ in this case is the same
as the parity in the previous case,
namely even.

(2)
We do some preliminary simplifications modulo $2$,
using \ref{rmktadmswitcht}:
$$
\eqalignno{
p(a',b',c') &=
\#K_1' \cdot (\pl(a',b') - 1)
+ \#K_2' \cdot (\#K_1' + \pl(a',c') - 1) \cr
&\hskip4em
+ \#K_3' \cdot (\#K_1' + \pl(b',c') + \#K_2' - 1) \cr
&\equiv
\#K_1' \cdot (\pl(a',a) + \pl(a,b) - 1)
+ \#K_2' \cdot (\#K_1' + \pl(a,c) - 1) \cr
&\hskip4em
+ \#K_3' \cdot (\#K_1' + \pl(b,c) + \pl(c,c') + \#K_2' - 1) \cr
&\equiv
\#K_1' \cdot \pl(a,b)
+ \#K_2' \cdot (\#K_1' + \pl(a,c) - 1) 
+ \#K_3' \cdot (\#K_1' + \pl(b,c) + \#K_2'). \cr
}
$$

Suppose that $i \in K_1$.
Then $a_i \not = b_i = A_i$,
and by \ref{lmtadm2},
$a_i \not = b_i = A_i = c_i$.
So $K_1' = K_1 \setminus \{i\}$.
Then $i \not \in K_2, K_2'$,
and exactly one of $K_3, K_3'$ contains $i$,
so that $\#K_3' = \#K_3 \pm 1$.
Thus
$$
\eqalignno{
p(a',b',c') &- p(a,b,c) + r =
\#K_1' \cdot \pl(a,b)
+ \#K_2' \cdot (\#K_1' + \pl(a,c) - 1) \cr
&\hskip4em
+ \#K_3' \cdot (\#K_1' + \pl(b,c) + \#K_2') \cr
&\hskip4em
- \#K_1 \cdot (\pl(a,b) - 1)
- \#K_2 \cdot (\#K_1 + \pl(a,c) - 1) \cr
&\hskip4em
- \#K_3 \cdot (\#K_1 + \pl(b,c) + \#K_2 - 1)
+ \ipl(a,c)
\cr
&\equiv
(\#K_1 + 1) \cdot \pl(a,b)
+ \#K_2 \cdot (\#K_1 + 1 + \pl(a,c) - 1) \cr
&\hskip4em
+ (\#K_3 + 1) \cdot (\#K_1 + 1 + \pl(b,c) + \#K_2) \cr
&\hskip4em
- \#K_1 \cdot (\pl(a,b) - 1)
- \#K_2 \cdot (\#K_1 + \pl(a,c) - 1) \cr
&\hskip4em
- \#K_3 \cdot (\#K_1 + \pl(b,c) + \#K_2 - 1)
+ \pl(a,c) - 1
\cr
&\equiv
\pl(a,b)
+ \#K_2 
+ \#K_1 + 1 + \pl(b,c) + \#K_2
- \#K_1
+ \pl(a,c) - 1,
\cr
}
$$
which is even.

So we may suppose that $i \not \in K_1$,
and since $p(a',b',c') - p(a,b,c) + r$ is even
if and only if $p(a,b,c) - p(a',b',c') + r$ is even,
by symmetry of the construction we may assume that $i \not \in K_1'$.
Then $a_i = b_i$ or $a_i = A_i$,
and $c_i = b_i$ or $c_i = A_i$.
If $a_i = b_i$,
then by \ref{lmtadm2},
$a_i = b_i \not = c_i = A_i$,
and necessarily $B_i = C_i = a_i$,
so that $i \in K_2$,
$K_2' = K_2 \setminus \{i\}$,
$K_3' = K_3$,
$K_1' = K_1$.
In this case,
$$
\eqalignno{
p(a',b',c') &- p(a,b,c) + r =
\#K_1' \cdot \pl(a,b)
+ \#K_2' \cdot (\#K_1' + \pl(a,c) - 1) \cr
&\hskip4em
+ \#K_3' \cdot (\#K_1' + \pl(b,c) + \#K_2') \cr
&\hskip4em
- \#K_1 \cdot (\pl(a,b) - 1)
- \#K_2 \cdot (\#K_1 + \pl(a,c) - 1) \cr
&\hskip4em
- \#K_3 \cdot (\#K_1 + \pl(b,c) + \#K_2 - 1)
+ \ipl(a,c)
\cr
&\equiv
\#K_1 \cdot \pl(a,b)
+ (\#K_2 - 1) \cdot (\#K_1 + \pl(a,c) - 1) \cr
&\hskip4em
+ \#K_3 \cdot (\#K_1 + \pl(b,c) + \#K_2 - 1) \cr
&\hskip4em
- \#K_1 \cdot (\pl(a,b) - 1)
- \#K_2 \cdot (\#K_1 + \pl(a,c) - 1) \cr
&\hskip4em
- \#K_3 \cdot (\#K_1 + \pl(b,c) + \#K_2 - 1)
+ \pl(a,c) - 1
\cr
&\equiv
- \#K_1 - \pl(a,c) + 1
- \#K_1
+ \pl(a,c) - 1,
\cr
}
$$
which is even.
So we may suppose that $a_i \not = b_i$,
so that $a_i = A_i$.
Since $c_i \not = a_i$,
by \ref{lmtadm2} then $A_i = a_i \not = b_i = c_i$.
It follows that $K_1' = K_1$,
$i \not \in K_2$,
$K_2' = K_2 \cup \{i\}$,
$K_3' = K_3$,
so that the
$p(a',b',c') - p(a,b,c) + r$ has the same parity as
the one in the case $a_i = b_i$,
which is even.


(3)
We do some preliminary simplifications modulo $2$,
using \ref{rmktadmswitcht}:
$$
\eqalignno{
p(a',b',c') &=
\#K_1' \cdot (\pl(a',b') - 1)
+ \#K_2' \cdot (\#K_1' + \pl(a',c') - 1) \cr
&\hskip4em
+ \#K_3' \cdot (\#K_1' + \pl(b',c') + \#K_2' - 1) \cr
&\equiv
\#K_1' \cdot (\pl(a,b) + \pl(b,b') - 1)
+ \#K_2' \cdot (\#K_1' + \pl(a,c) + \pl(c,c') - 1) \cr
&\hskip4em
+ \#K_3' \cdot (\#K_1' + \pl(b,c) + \#K_2' - 1) \cr
&\equiv
\#K_1' \cdot \pl(a,b)
+ \#K_2' \cdot (\#K_1' + \pl(a,c)) \cr
&\hskip4em
+ \#K_3' \cdot (\#K_1' + \pl(b,c) + \#K_2' - 1). \cr
}
$$

If $i \in K_1$,
then $a_i \not = b_i = A_i$.
Since $b_i \not = c_i$,
by \ref{lmtadm2} then $c_i = a_i \not = b_i = A_i$,
so that $K_1' = K_1 \setminus \{i\}$,
$i \not \in K_2, K_3$,
$K_2' = K_2 \cup \{i\}$ and $K_3' = K_3$.
Thus here
$$
\eqalignno{
p(a',b',c') &- p(a,b,c) + r =
\#K_1' \cdot \pl(a,b)
+ \#K_2' \cdot (\#K_1' + \pl(a,c)) \cr
&\hskip4em
+ \#K_3' \cdot (\#K_1' + \pl(b,c) + \#K_2' - 1)
\cr
&\hskip4em
- \#K_1 \cdot (\pl(a,b) - 1)
- \#K_2 \cdot (\#K_1 + \pl(a,c) - 1) \cr
&\hskip4em
- \#K_3 \cdot (\#K_1 + \pl(b,c) + \#K_2 - 1)
+ \ipl(b,c)
\cr
&\equiv
(\#K_1 - 1) \cdot \pl(a,b)
+ (\#K_2 + 1) \cdot (\#K_1 - 1 + \pl(a,c)) \cr
&\hskip4em
+ \#K_3 \cdot (\#K_1 - 1 + \pl(b,c) + \#K_2 + 1 - 1)
\cr
&\hskip4em
- \#K_1 \cdot (\pl(a,b) - 1)
- \#K_2 \cdot (\#K_1 + \pl(a,c) - 1) \cr
&\hskip4em
- \#K_3 \cdot (\#K_1 + \pl(b,c) + \#K_2 - 1)
+ \pl(b,c) - 1
\cr
&\equiv
- \pl(a,b)
+ \#K_1 - 1 + \pl(a,c))
- \#K_1
+ \pl(b,c) - 1,
\cr
}
$$
which is even.

So we may assume that $i \not \in K_1$,
and since $p(a',b',c') - p(a,b,c) + r$ is even
if and only if $p(a,b,c) - p(a',b',c') + r$ is even,
by symmetry of the construction we may assume that $i \not \in K_1'$.
Then either $a_i = b_i$ or $a_i = A_i$,
and either $a_i = c_i$ or $a_i = A_i$.
If $a_i = A_i$,
then $i \not \in K_2', K_2$,
$K_1' = K_1$,
$K_2' = K_2$,
and exactly one of $K_3, K_3'$ has $i$,
so that $\#K_3' = \#K_3 \pm 1$.
Thus modulo $2$:
$$
\eqalignno{
p(a',b',c') &- p(a,b,c) + r =
\#K_1' \cdot \pl(a,b)
+ \#K_2' \cdot (\#K_1' + \pl(a,c)) \cr
&\hskip4em
+ \#K_3' \cdot (\#K_1' + \pl(b,c) + \#K_2' - 1)
\cr
&\hskip4em
- \#K_1 \cdot (\pl(a,b) - 1)
- \#K_2 \cdot (\#K_1 + \pl(a,c) - 1) \cr
&\hskip4em
- \#K_3 \cdot (\#K_1 + \pl(b,c) + \#K_2 - 1)
+ \ipl(b,c)
\cr
&\equiv
\#K_1 \cdot \pl(a,b)
+ \#K_2 \cdot (\#K_1 + \pl(a,c)) \cr
&\hskip4em
+ (\#K_3 + 1) \cdot (\#K_1 + \pl(b,c) + \#K_2 - 1)
\cr
&\hskip4em
- \#K_1 \cdot (\pl(a,b) - 1)
- \#K_2 \cdot (\#K_1 + \pl(a,c) - 1) \cr
&\hskip4em
- \#K_3 \cdot (\#K_1 + \pl(b,c) + \#K_2 - 1)
+ \pl(b,c) - 1
\cr
&\equiv
\#K_1 + \pl(b,c) + \#K_2 - 1
+ \#K_1
+ \#K_2
+ \pl(b,c) - 1,
\cr
}
$$
which is even.
So we may assume that $a_i \not = A_i$.
Then $a_i = b_i = c_i$,
which contradicts the assumption that $b_i \not = c_i$.
\qed

\prop[propGBred]
Let $S$ be a $t$-signed set,
and let $a,b,c \in S$ be mutually connected.
Then all reductions of $x_a x_b x_c$
in the lexicographic order with respect to $\CG{t}{S}$
reduce to the same term,
and the term is of the form $(-1)^u x_A x_B x_C$ for some integer $u$,
unique up to parity.
\endb

\proof
Set
$G' = \{f_{K,d,e} : d, e \in S, K \subseteq \{i \in [t]: d_i \not = e_i\}\}$.
This set is reminiscent of $\CG{t}{S}$,
the only difference is that each binomial in $G'$
has the two coefficients $1, -1$,
whereas some binomials in $\CG t S$ have coefficients $1, 1$
and others have $1 , -1$.
By Theorem~4.5 of \cite{ST},
$x_a x_b x_c$ and all of its reductions with respect to $G'$
reduce to some minimal monomial $x_A x_B x_C$
with respect to $G'$.
By the form of $G'$ and $\CG{t}{S}$,
then with respect to $\CG{t}{S}$,
$x_a x_b x_c$ reduces to $(-1)^{u} x_A x_B x_C$
for some integer $u$,
and we need to show that up to parity $u$ is uniquely determined
regardless of what the reduction steps are.
In fact,
we show that $u$ has the same parity as
$p(a,b,c)$ from \ref{lmGB}.
If $a = A$, $b = B$, and $c = C$,
then $p(a,b,c) = 0$
and $x_a x_b x_c$ is in the reduced form and it does not reduce any further,
so the conclusion holds trivially.

Now suppose that in the first step of the reduction,
we reduce $x_a x_b x_c$ with respect to $h_{L,a,b}$
for some non-empty $L \subseteq \{i \in [t]: a_i \not = b_i\}$.
Then $x_a x_b x_c$ is reduced to
$(-1)^{\#L \cdot \ipl(a,b)} x_{\switch(L,a,b)} x_{\switch(L,b,a)} x_c$.
This is a proper reduction,
so by induction on the order (in the lexicographic order),
$x_{\switch(L,a,b)} x_{\switch(L,a,b)} x_c$
reduces to $(-1)^{p(\switch(L,a,b), \switch(L,b,e), c)} x_A x_B x_C$.
Hence via this reduction
we have $u =
p(\switch(L,a,b), \switch(L,b,e), c)
+ \#L \cdot \ipl(a,b)$.
Write $L = \{l_1, \ldots, l_k\}$.
Then we have the following modulo $2$:
$$
\eqalignno{
u - p(a,b,c)
&=
p(\switch(L,a,b), \switch(L,b,a), c) - p(a,b,c)
+ \#L \cdot \ipl(a,b) \cr
&\equiv
\sum_{i=1}^k \biggl(
p(\switch(\{l_1, \ldots, l_i\},a,b),
\switch(\{l_1, \ldots, l_i\},b,a),
c) \cr
&\hskip5em
- p(\switch(\{l_1, \ldots, l_{i-1}\},a,b),
\switch(\{l_1, \ldots, l_{i-1}\},b,a),
c)
\cr
&\hskip5em
+ 
\ipl(l_i,\switch(\{l_1, \ldots, l_{i-1}\},a,b),
\switch(\{l_1, \ldots, l_{i-1}\},b,a))
\biggr), \cr
}
$$
which is even by \ref{lmGBcont}.
A very similar proof shows the same conclusion
if we first reduce $x_a x_b x_c$ with respect to $h_{L,a,c}$
or $h_{L,b,c}$.
\qed

%


\thm[thmGB]
Let $S$ be a $t$-signed set.
Then the set $\CG{t}{S}$
is a (non-minimal) Gr\"obner basis for $\CJtS{t}{S}$ in the lexicographic order.
\endb

\proof
%
A note:
``S" in ``S-polynomial" below is unrelated to the $t$-signed set $S$.

Since $h_{K,a,b}
= x_a x_b
- (-1)^{\#K \cdot \ipl(a,b)} x_{\switch(K,a,b)} x_{\switch(K,b,a)}
= \pm h_{K, \switch(K,a,b), \switch(K,b,a)}$,
the ideal generated by $\CG{t}{S}$
is the same as the ideal generated by
$G = \{h_{K,a,b} \in \CG{t}{S}:
x_a x_b > x_{\switch(K,a,b)} x_{\switch(K,b,a)}\}$.
It suffices to prove that $G$ is a Gr\"obner basis.

We have to prove that for any $h_{K,a,b}, h_{L,c,d} \in G$,
the S-polynomial $S( h_{K,a,b}, h_{L,c,d})$
reduces to $0$ in the Gr\"obner basis sense with respect to $G$.

If $x_a x_b$ and $x_c x_d$ have no variables in common,
then by standard facts about Gr\"obner bases
their S-polynomial reduces to $0$.
So we may assume that $a = d$.
Then $a, b, c, d$ are all connected in $S$.

Suppose first that in addition $b = c$.
Then
$S(h_{K,a,b}, h_{L,b,a})$ equals
$$
-(-1)^{\#K \cdot \ipl(a,b)} x_{\switch(K,a,b)} x_{\switch(K,b,a)}
+ (-1)^{\#L \cdot \ipl(a,b)} x_{\switch(L,a,b)} x_{\switch(L,b,a)}.
$$
If this is $0$,
we are done,
otherwise,
this equals
$$
\eqalignno{
-(-1)^{\#K \cdot \ipl(a,b)}
&\left(
x_{\switch(K,a,b)} x_{\switch(K,b,a)}
+ (-1)^{(\#L - \#K) \cdot \ipl(a,b)} x_{\switch(L,a,b)} x_{\switch(L,b,a)}
\right) \cr
&= \pm
\bigl(
x_{\switch(K,a,b)} x_{\switch(K,b,a)}
\cr
&\hskip2em
+ (-1)^{(\#(L \setminus K) + \#(K \setminus L)) \cdot \ipl(a,b)}
x_{\switch((L \setminus K) \cup (K \setminus L),\switch(K,a,b),\switch(K,b,a))}
\cr
&\hskip2em
\hphantom{+ (-1)^{(\#(L \setminus K) + \#(K \setminus L)) \cdot \ipl(a,b)}}
\cdot
x_{\switch((L \setminus K) \cup (K \setminus L),\switch(K,b,a),\switch(K,a,b))}
\bigr) \cr
&= \pm
h_{(L \setminus K) \cup (K \setminus L),\switch(K,a,b),\switch(K,b,a))}, \cr
}
$$
which is a scalar multiple of an element in $G$,

So we may assume that $a = d$ and $b \not = c$.
Then
$$
\eqalignno{
S(&h_{K,a,b}, h_{L,a,c}) \cr
&=
-(-1)^{\#K \cdot \ipl(a,b)} x_{\switch(K,a,b)} x_{\switch(K,b,a)} x_c
+ (-1)^{\#L \cdot \ipl(a,c)} x_{\switch(L,a,c)} x_{\switch(L,c,a)} x_b. \cr
}
$$
Both
$(-1)^{\#K \cdot \ipl(a,b)} x_{\switch(K,a,b)} x_{\switch(K,b,a)} x_c$
and $(-1)^{\#L \cdot \ipl(a,c)} x_{\switch(L,a,c)} x_{\switch(L,c,a)} x_b$
are reductions of $x_a x_b x_c$,
so that by \ref{propGBred},
S-polynomial reduces to $0$ with respect to $\CG{t}{S}$.
\qed

\thm[thmGBprime]
If $S$ is a $t$-signed set,
then the ideals $\CJtS{t}{S}$ and $\CQ{t}{S}$ are prime.
\endb

\proof
This proof mimics that of Theorem~6.3 in~\cite{ST}.
By faithfully flatness
we may assume without loss of generality that
the underlying field is algebraically closed.

Let $S = S_1 \cup \cdots \cup S_k$
be a partition of $S$ into equivalence classes
with respect to connectedness.
Each $S_i$ is $t$-signed.
Then $\CJtS{t}{S} = \cup_i \CJtS{t}{S_i}$
and $\CQ{t}{S} = \sum_i \CJtS{t}{S_i} + \CVar{t}{S}$,
and the generators of
$\CJtS{t}{S_1}$,
$\ldots$,
$\CJtS{t}{S_k}$,
and $\CVar{t}{S}$
use disjoint variables.
By a well-known fact,
it suffices to prove that each $\CJtS{t}{S_i}$ is a prime ideal.
By renaming we now assume that $t$-signed $S = S_i$
is one equivalence class under connectedness.

Here we quote Lemma 6.2 from~\cite{ST}:
Suppose that $a_1, \ldots, a_r, b_1, \ldots, b_r \in S$
have the property that for all $i = 1, \ldots, t$,
up to order, the multiset $\{a_{1i}, a_{2i}, \ldots, a_{ri}\}$
is the same as the multiset $\{b_{1i}, b_{2i}, \ldots, b_{ri}\}$,
and such that, up to order, the multiset
$\{(a_{1,t+1}, a_{1,t+2}, \ldots, a_{1,n}),
(a_{2,t+1}, a_{2,t+2}, \ldots, a_{2,n}),
\ldots,
(a_{r,t+1}, a_{r,t+2}, \ldots, a_{r,n})\}$
is the same as the multiset
$\{(b_{1,t+1}, b_{1,t+2}, \ldots, b_{1,n}),
(b_{2,t+1}, b_{2,t+2}, \ldots, b_{2,n}),
\ldots,
(b_{r,t+1}, b_{r,t+2}, \ldots, b_{r,n})\}$.
Then in the lexicographic order,
$x_{a_1} x_{a_2} \cdots x_{a_r} - x_{b_1} x_{b_2} \cdots x_{b_r}$
reduces with respect to $\{f_{K,a,b} : K \subseteq [t], a, b \in S\}$
to $0$.

A consequence is that under the conditions in the previous paragraph,
either
$x_{a_1} \cdots x_{a_r} - x_{b_1} \cdots x_{b_r}$
or $x_{a_1} \cdots x_{a_r} + x_{b_1} \cdots x_{b_r}$
reduces to $0$ with respect to $G$.

Suppose that $\CJtS{t}{S}$ is not a prime ideal.
Since this is a binomial ideal,
by Eisenbud--Sturmfels~\cite{ESbinom}
there exists a zerodivisor modulo $\CJtS{t}{S}$ of the form $\alpha - c \beta$
for some monomials $\alpha$ and $\beta$
and some possibly zero coefficient $c$ in the base field.
Then there exists $f$ not in $\CJtS{t}{S}$ such that $(\alpha - c \beta) \cdot f \in \CJtS{t}{S}$.
Without loss of generality we may assume that $\alpha$, $\beta$,
and each monomial in $f$ is reduced with respect to~$\CG{t}{S}$.
Write $f = c_1 m_1 + c_2 m_2 + \cdots + c_k m_k$,
where the $c_i$ are non-zero elements of the underlying field
and the $m_i$ are monomials.
We may assume that $m_1 > m_i$ for all $i$
and $\alpha > \beta$ in the lexicographic order.
Since $(\alpha - c \beta) \cdot f$ reduces to $0$
with respect to~$\CG{t}{S}$,
$\alpha m_1$ must reduce with respect to $\CG{t}{S}$ up to sign
to the same monomial as some other monomial in $(\alpha - c \beta) \cdot f$.
If it reduces to the same monomial as $\alpha m_i$ for some $i > 1$,
then by the previous paragraph,
$m_1$ and $m_i$ reduce up to sign to the same monomial,
which contradicts the assumption on the monomials in $f$ being reduced.
Hence $\alpha m_1$
reduces with respect to $\CG{t}{S}$ to the same monomial as some $\beta m_{i_2}$,
and by the reduced assumption on $\alpha - c \beta$ necessarily $i_2 > 1$
and $c$ is not zero.
Similarly, $\alpha m_{i_2}$ reduces to the same monomial as $\beta m_{i_3}$,
for some $i_3 \not = i_2$,
and more generally for all $s$,
$\alpha m_{i_s}$ reduces to the same monomial as $\beta m_{i_{s+1}}$
for some $i_s \not = i_{s+1}$.
Necessarily some $i_s$ must equal some $i_j$ for $s < j$.
Hence $\alpha^{j-s+1} m_{i_s} \cdots m_{i_j}$ reduces to the same monomial as
$\beta^{j-s+1} m_{i_s} \cdots m_{i_j}$,
whence again by the previous paragraph,
$\alpha^{j-s+1}$ reduces to the same monomial as $\beta^{j-s+1}$,
and even $\alpha$ reduces to the same monomial as $\beta$,
which is a contradiction.
Thus $\CJtS{t}{S}$ is a prime ideal.
\qed

In particular,
$\CJtS{t}{N}$ is prime
in case $N = [2] \times [2] \times [2]$.
We give here an easier proof in this special case based on results of \cite{ST},
but this easier proof does not generalize to arbitrary~$N$.
Namely,
$\CJtS{t}{N}$
$= \left(f_{i,a,b} : d(a,b) = 3, i \in [t]\right)$
$+ \left(g_{i,a,b} : d(a,b) = 2, i \in [t], a_i \not = b_i\right)$.
Let
$\varphi : R \rightarrow R$
be the ring isomorphism that restricts to the identity map
on the underlying field
and maps $x_{ijk}$
to $-x_{ijk}$ if $i=j=k$
and to $x_{ijk}$ otherwise.
It is straightforward to see that $\varphi$
takes $\CJtS{t}N$ onto $\CItS{t}N$.
But by~\cite{ST}, $\CItS{t}N$ is a prime ideal (for all sizes of~$N$),
and since isomorphisms map prime ideals to prime ideals,
the conclusion follows.

\section{sectminprimescont}{{Prime ideals minimal over $\CJtS{t}{S}$}}

We prove in this section
that the maximal $t$-signed sets correspond precisely
to prime ideals minimal over $\CJ{t}$.
In~\cite{ST} it was proved
that the maximal $t$-switchable sets correspond precisely
to prime ideals minimal over $\CI{t}$.
The flow of the proofs resembles those in~\cite{ST},
but again,
here in addition parity has to be checked and controlled.

\thm[thmCJminprimesone]
If $P$ is a prime ideal minimal over $\CJ{t}$,
then $P = \CQ{t}{S}$ for some maximal $t$-signed set $S$.
\endb

\proof
Let $S$ be the set of all $a \in [r_1] \times \cdots \times [r_n]$
such that $x_a \not \in P$.

Let $a, b \in S$ have $d(a,b) = 2$ and $a_i\neq b_i$ for some $i\in[t]$.
Since $P$ contains $\CJ{t}$ and $i\in[t]$,
$P$ contains $h_{S,i,a,b}
= x_a x_b - (-1)^{\ipl_S(a,b)}  x_{\switch(i,b,a)} x_{\switch(i,a,b)}$.
Since $a, b \in S$,
then $x_a x_b \not \in P$,
so that necessarily
$x_{\switch(i,b,a)} x_{\switch(i,a,b)} \not \in P$,
and hence
$\switch(i,b,a), \switch(i,a,b) \in S$.
This proves that $S$ is $t$-switchable.

We next prove that $S$ is $t$-signed.
It suffices to prove that every connected component $S_0$ of $S$
is $t$-signed.
If all elements of $S_0$ share the same first $t$ coordinates,
or if all elements in $S_0$ only differ in the $i$th component
for some $i \in [t]$,
then $S_0$ is $t$-signed.
So we may suppose that the distance between some elements $a, b$
of $S_0$ is at least $2$
and that for some $i \in [t]$,
$\{c_i: c\in S_0\}$ has more than one element.
By \ref{rmktadmswitcht},
by possibly replacing $a$ by $\switch(i,a,c)$ for some $c \in S_0$,
we may assume that $a$ and $b$ differ in the $i$th component.
Now let $c, e \in S_0$ be arbitrary.
If there are paths from $c$ to $e$ whose lengths have different parities,
then there are paths from $a$ to $b$ via paths from $c$ to $e$
whose lengths have different parities.
Then by \ref{proppathparity},
$\CJ{t}$ and hence $P$ contain monomials
whose variables have subscripts in $S$.
But $P$ is a prime ideal,
so $P$ contains a factor $x_c$ for some $c \in S$,
which is a contradiction.
This proves that $S$ is $t$-signed.

By \ref{propCJtinCQ},
$\CJ{t} \subseteq \CQ{t}{S}$,
and by \ref{thmGBprime},
$\CQ{t}{S}$ is a prime ideal.

We next prove that $\CQ{t}{S} \subseteq~P$.
By the definition of $S$, $\CVar{t}{S} \subseteq P$.
Let $h_{S,i,a,b} \in \CJtS{t}{S}$,
with $i\in[t]$ and $a$ and $b$ connected in $S$
such that $a_i \not = b_i$.
If $d(a,b) = 1$,
then $h_{S,i,a,b} = 0$,
so it is an element of $P$.
So we may assume that $d(a,b) \ge 2$.
By the definition of connectedness,
there exist elements $c_0 = a, c_1, \ldots, c_k, c_{k+1}, b \in S$
such that for all $j = 1, \ldots, k$,
$c_{j-1}$ and $c_j$ differ only in one position.
Without loss of generality $d(c_k,b) = 2$.
After omitting any repetitions in
$c_0, \switch(i,c_1,a), \ldots, \switch(i,c_{k+1},a)$,
all the $i$th components on the list equal $a_i$.
By \ref{rmktadmswitcht},
this list is in $S$.
Note that $2 \le d(\switch(i,c_k,a),b) \le 3$.
If $d(\switch(i,c_k,a),b) = 3$,
then we take the list
$c_0, \switch(i,c_1,a), \ldots, \switch(i,c_{k+1},a)$
with redundancies removed,
and the last element on the list differs from $b$ in $2$ entries,
and if $d(\switch(i,c_k,a),b) = 2$,
then we take the list
$c_0, \switch(i,c_1,a), \ldots, \switch(i,c_k,a)$
with redundancies removed.
In either case,
after renaming we have a path
$a = c_0, c_1, \ldots, c_k$
with all $i$th components being $a_i$ and $d(c_k,b) = 2$.
Then by~\ref{lmvarsmult},
$x_{c_1} \cdots x_{c_k} h_{S,i,a,b} \in \CJ{t} \subseteq P$,
and since $x_{c_j} \not \in P$,
it follows that $h_{S,i,a,b} \in P$,
as desired.
Thus $\CJtS{t}{S} \subseteq \CQ{t}{S} \subseteq P$.
Since $\CQ{t}{S}$ is a prime ideal,
by minimality of $P$, $\CQ{t}{S} = P$.

Finally,
let $T$ be $t$-signed and properly containining $S$.
Then $\CVar{t}{T} \subsetneq \CVar{t}{S}$,
and so $\CQ{t}{T} \not = \CQ{t}{S}$.
By~\ref{propCJtinCQ},
$\CQ{t}{T}$ contains $\CJ{t}$,
and by \ref{thmGBprime},
$\CQ{t}{T}$ is a prime ideal.
This combined with the fact that $P = \CQ{t}{S}$ is minimal over $\CJ{t}$,
implies that $\CQ{t}{T} \not\subseteq \CQ{t}{S}$.
Therefore $\CQ{t}{S}$ and $\CQ{t}{T}$ are incomparable.
Thus $S$ is a maximal $t$-signed set.

\thm[thmCJminprimes]
The set of prime ideals minimal over~$\CJ{t}$
equals the set of ideals of the form $\CQ{t}{S}$
as $S$ varies over maximal $t$-signed sets.
\endb

\proof
Let $S$ be a maximal $t$-signed set.
Then by \ref{propCJtinCQ} and \ref{thmGBprime},
$\CQ{t}{S}$ is a prime ideal that contains~$\CJ{t}$.
Suppose that $P$ is a prime ideal that is minimal over $\CJ{t}$
and is contained in $\CQ{t}{S}$.
By~\ref{thmCJminprimesone},
$P = \CQ{t}{T}$ for some maximal $t$-signed set~$T$.
Since $\CQ{t}{T} \subseteq \CQ{t}{S}$,
necessarily $\CVar{t}{T} \subseteq \CVar{t}{S}$,
so that $S \subseteq T$.
But then maximality of $S$ forces $\CQ t T = \CQ t S$.
Thus $\CQ{t}{S}$ is a minimal prime ideal of~$\CJ{t}$.
\ref{thmCJminprimesone} proves the other direction.
\qed

\section{sectadmanalysis}{Structure of $t$-signed sets
and  the radical of $\CJ t$}

The aim of this section is to establish the structure of $t$-signed sets
in more detail,
with the bigger goal of then determining the radical of $\CJ t$.
However,
our description of this radical is only indirect.

\lemma[lm2-box]
For $i = 1, \ldots, n$,
let $S_i \subseteq [r_i]$.
Suppose that one of the following conditions holds:
\ritem
$|S_1| = \cdots = |S_t| = 1$.
\ritem
There exists $i$ such that for all $j \not = i$,
$|S_j| = 1$.
\ritem
For all $i$,
$|S_i| \le 2$.

\noindent
Then $S = S_1 \times \cdots \times S_n$ is a $t$-signed set
consisting of one equivalence class only.
\endb

\proof
It is clear that $S$ forms one equivalence class.
Conditions (1) and (2) above correspond precisely to
conditions (1) and (2) in \ref{deftadm}.

Now assume that condition (3) above holds.
If $a, b \in S$ such that $d(a,b) = 2$
and $i \in [t]$ with $a_i \not = b_i$,
then $\switch(i,a,b), \switch(i,b,a) \in S$ as well,
so that $S$ is $t$-switchable.
Any path from an arbitrary $a$ to an arbitrary $b$ in $S$
may involve several switches in each of the entries,
and the parity of the number of switches  in the $k$th entry
is $1$ if $a_k \not = b_k$ and $0$ otherwise.
So the parity of each path,
namely the parity of the sum of the switches in all the components,
is uniquely determined.
Thus $S$ is $t$-signed.
\qed

By \ref{ex223-1adm}~(4),
the $t$-switchable sets
can have a form different from the forms given in \ref{lm2-box} above.


\lemma[lmnot2by3]
Let $U$ be a subset of $N$
that is contained in a $2 \times 3$ submatrix of $N$
with entries varying in the $i$ and $j$th components only.
Suppose that $U$ is not contained in a $2 \times 2$ submatrix
or in a $1 \times 3$ submatrix,
and that either $i$ or $j$ is in $[t]$.
Then $U$ is not a subset of any $t$-signed set.

Conversely,
if $U$ is not a subset of any $t$-signed $S$,
then the smallest $t$-switchable set containing $U$ contains
a $2 \times 3$ submatrix of $N$
with entries varying in the $i$ and $j$th components
with either $i$ or $j$ in $[t]$.
\endb

\proof
Suppose that $U$ is a subset of a $t$-signed set $S$.
Then since $S$ is $t$-switchable,
by \ref{rmktadmswitcht},
the whole $2 \times 3$ matrix is in $S$,
and even all its entries are in the same connected equivalence class $S_0$.
So,
let $a, b, c \in S_0$ be in a $1 \times 3$ submatrix.
Then $a, b, c$ and $a, c$ are both paths from $a$ to $b$ in $S_0$,
but then $S$ is not $t$-signed,
which is a contradiction.
The proof of the converse is similar.
\qed

\defin[defmonom]
For any finite multiset $M$ of elements of $N$,
define $x_M = \prod_{a \in M} x_a$.
(When $M$ is a multiset,
$x_M$ allows for products of powers of variables,
but when $M$ is a set,
then $x_M$ is square-free.)
\endb

\lemma[lmMinrad]
Let $\listmonoms$ be the set of all subsets of $N$
that are not contained in any $t$-signed sets.
A monomial $x_M$ is in $\sqrt{\CJ t}$
if and only if $M \in \listmonoms$.
\endb

\proof
By \ref{thmCJminprimes},
$x_M \in \sqrt{\CJ t}$
if and only if for all maximal $t$-signed $S$,
$x_M \in \CQ t S$,
which by the structure of these prime ideals holds if and only if
for all $S$,
$M$ is not a subset of $S$.
This is the same as saying that $M \in \listmonoms$.
%
\qed


\cor[cormonomrad]
Let $\listmonoms_i$ be the set of all sets in $\listmonoms$ with $i$ elements.
Then $(x_M: M \in \listmonoms) + \CJ t$
$= (x_M: M \in \listmonoms_3) + \CJ t$.
\endb

\proof
For any $a, b \in N$,
by \ref{lm2-box},
$\{a_1,b_1\} \times \cdots \times \{a_n,b_n\}$ is $t$-signed,
so that $\listmonoms_0 = \listmonoms_1 = \listmonoms_2 = \emptyset$.
It suffices to prove that
$(x_M: M \in \listmonoms) \subseteq (x_M: M \in \listmonoms_3) + \CJ t$.
Let $M \in \listmonoms$.
By \ref{lmnot2by3},
the smallest $t$-switchable set containing $M$ contains a $2 \times 3$
submatrix with fixed coordinates $i, j$
such that $\{i,j\} \cap [t] \not = \emptyset$.
It need not be the case that this $2 \times 3$ submatrix
contains $3$ elements of $M$
that do not lie in a $1 \times 3$ or $2 \times 2$ submatrix.
However,
let $S$ be the smallest subset of $N$
containing $S$
and such that whenever $a, b \in S$ with $d(a,b) = 2$,
then $\switch(i,a,b) \in S$.
The assumption is that elements of $S$
fill that $2 \times 3$ submatrix.
But generating $\switch(i,a,b), \switch(i,b,a)$ from $a, b$
is on the algebraic side
the same as subtracting multiples of generators of $\CJ t$,
which says that after repeatedly subtracting from $x_M$
specific elements of $\CJ t$,
we get a monomial $x_{M'}$
that has a factor
$x_a x_b x_c$
such that the smallest $t$-switchable set containing $a, b, c$
contains that $2 \times 3$ submatrix.
This proves the corollary.
\qed

It is not true that $(x_M: M \in \listmonoms) = (x_M: M \in \listmonoms_3)$.
For example,
if $n = t = 3$,
then $M = \{(1,1,1), (2,1,1), (3,1,1), (1,2,2)\}$ is in $\listmonoms$.
By \ref{lm2-box},
the proper subsets
$\{(1,1,1), (3,1,1), (1,2,2)\}$
$\{(1,1,1), (2,1,1), (1,2,2)\}$
$\{(1,1,1), (2,1,1), (3,1,1)\}$
are not in~$\listmonoms$,
and $\{(2,1,1), (3,1,1), (1,2,2)\}$ is 
$t$-signed with exactly two equivalence classes under connectedness
($\{(2,1,1), (3,1,1)\}$ and $\{(1,2,2)\}$),
and so this last subset of $M$ is also not in $\listmonoms$.
Thus $x_M \not \in (x_{M'}: M' \in \listmonoms_3)$.

\thm
Let $\pairslistvarsTswrad$ be the set of all pairs $(M,M')$
of finite lists of elements of $N$
for which there exists an integer $v(M,M')$
such that for all $t$-signed $S$,
$x_M - (-1)^{v(M,M')} x_{M'} \in \CQ t S$.
The radical of $\CJ t$ is generated by all elements
$x_M - (-1)^{v(M,M')} x_{M'}$ for $(M, M') \in \pairslistvarsTswrad$,
and by all $x_M$ for $M \in \listmonoms_3$.
\endb

\proof
Let $J' = (x_M - (-1)^{v(M,M')} x_{M'}: (M, M') \in \pairslistvarsTswrad)$
and $J'' = (x_M : M \in \listmonoms_3)$.
Clearly $\CJ t \subseteq J'$.

By \ref{thmCJminprimes},
$J' \subseteq \sqrt{\CJ t}$,
and
by \ref{lmMinrad},
$J'' \subseteq \sqrt{\CJ t}$.
Thus $J' + J'' \subseteq \sqrt{\CJ t}$.

For the other inclusion,
first  assume that the underlying field is algebraically closed.
By~\cite{ESbinom},
$\sqrt{\CJ t}$ is generated by binomials and monomials.
Let $f = x_M - c x_{M'} \in \sqrt{\CJ t}$
for some finite lists $M, M'$ and some possibly zero scalar $c$.

Suppose that $x_M \in \sqrt{\CJ t}$.
Then by \ref{lmMinrad},
the set $M$ is in $\listmonoms$,
so that by \ref{cormonomrad},
$x_M \in \CJ t + J''$,
and either $c = 0$ or similarly $x_{M'} \in \CJ t + J''$,
whence $f \in \CJ t + J''$.

Next assume that $x_M \not \in \sqrt{\CJ t}$.
Thus $c$ is non-zero.
By \ref{thmCJminprimes},
$f \in \CQ t S$ for all $t$-signed $S$,
and for at least one such $S$,
$x_M \not \in \CQ t S$.
By the structure of these prime ideals,
$x_M - (-1)^{p_S} x_{M'} \in \CQ t S$
for some integer $p_S$ depending on $S$.
But also $x_M - c x_{M'} \in \sqrt{\CJ t} \subseteq \CQ t S$,
so that $c = (-1)^{p_S}$.
Thus by assumption for any $t$-signed set $T$,
$x_M - (-1)^{p_S} x_{M'} \in \sqrt{\CJ t} \subseteq \CQ t T$,
hence $f \in J'$.

This proves that $\sqrt{\CJ t} = J' + J''$
whenever the underlying field is algebraically closed.

For an arbitrary underlying field,
let $R'$ be $R$ tensored with the algebraic closure of the field.
Then by above,
 $\sqrt{\CJ t R'} = (J' + J'') R'$,
which has generators in $R$,
whence by faithful flatness of $R'$ over $R$,
also $\sqrt{\CJ t} = J' + J''$.
\qed

If $N = [r_1] \times [r_2]$,
by~\cite{LS},
$\sqrt{\CJ t}
= \CJ t
+ (x_M: M \in \listmonoms_3)$.
However,
for $N = [2] \times [2] \times [2] \times [2]$,
$\sqrt{\CJ t}$ properly contains
$\CJ t + (x_M: M \in \listmonoms_3)$.

\section{sectjuliahat}{Related ideals I:
generated by $g_{i,a,b}$ when $d(a,b) = 3$}


In this section we introduce another notion of distance
$$
d_t(a,b) = \#\{i \in [t]: a_i \not = b_i\},
$$
and the set
$\{i \in [t]: a_i \not = b_i\}$ will be denoted as $D_t(a,b)$.

\defin[defJhat]
Let $t \le n$
and let $\ChJ{t}$ be the ideal generated by $g_{i,a,b}$
with $a, b \in N$
such that $d(a,b) = d_t(a,b) = 3$ and $i \in D_t(a,b)$.
\endb

\lemma[lemmamonogens]
$\ChJ{t}$ is generated by monomials $x_a x_b$
such that $d(a,b) = d_t(a,b)=3$ and $i \in D_t(a,b)$.
Thus $\ChJ{t}$ is a radical ideal.
\endb

\proof
If $t \le 2$,
$\ChJ{t}$ is the zero ideal,
and so is the ideal generated by the non-empty set of specified monomials.
So we may assume that $t \ge 3$.

Let $a, b \in N$ satisfy $d(a,b) = d_t(a,b) = 3$.
Let $D_t(a,b) = \{i,j,k\}$.
Then $g_{i,a,b}, g_{j,a,b}, g_{k,a,b}$ are among the generators of $\ChJ{t}$,
but $g_{i,a,b} - g_{j,a,b}
= x_{\switch(i,a,b)} x_{\switch(i,b,a)}
- x_{\switch(j,a,b)} x_{\switch(j,b,a)}$,
and
$g_{k,a,b}
= x_{\switch(i,a,b)} x_{\switch(i,b,a)}
+ x_{\switch(j,a,b)} x_{\switch(j,b,a)}$,
so that as the characteristic of the underlying field is not $2$,
$\ChJ{t}$ contains
$x_{\switch(i,a,b)} x_{\switch(i,b,a)}$
and $x_{\switch(j,a,b)} x_{\switch(j,b,a)}$,
whence it also contains $x_a x_b$.
This proves one inclusion,
and the other is trivial.
So $\ChJ t$ is generated by square-free monomials,
so it is a radical ideal.
\qed

The following is now immediate:

\thm[thmminprimeshat]
Let ${\bf \widehat S}$ be the collection of all subsets $S$ of $N$
such that for any $a,b \in N$ with $d(a,b) = d_t(a,b) = 3$,
at least one of $a$ and $b$ is not in $S$.
Then the prime ideals that are minimal over $\ChJ{t}$
are the minimal ideals in $\{\CVar{t}{S}: S \in {\bf \widehat S}\}$
(see \ref{defadmPQ}).
\qed
\endb


It is clear that when $t = 3$,
the set ${\bf \widehat S}$
consists of sets $S$ containing points no two of which have distance $3$.
In order to get the minimal primes,
we need $\CVar{3}{S}$ minimal possible,
so we want $S$ maximal possible.

%

\example
If $N = [r_1] \times [r_2] \times [r_3]$ and $t = 3$,
the minimal prime ideals over $\ChJ{3}$
correspond to the sets $S$,
whose geometric descriptions are as followis:
\ritem
For $i = 1, 2, 3$,
$S$ consists of all elements with a fixed $i$th coordinate;
there are $r_1 + r_2 + r_3$ such sets.
\ritem
$S$ consists of four elements in a $2 \times 2 \times 2$ subhypermatrix
consisting of one of the corners plus its adjacent neighbors;
there are ${r_1 \choose 2} {r_2 \choose 2} {r_3 \choose 2}$
such subhypermatrices,
and each one has eight such sets.
\ritem
$S$ consists of four elements in a $2 \times 2 \times 2$ subhypermatrix
all of whose pairs have distance $2$;
there are ${r_1 \choose 2} {r_2 \choose 2} {r_3 \choose 2}$
such subhypermatrices,
and each one has two such sets.

Thus in total there are
$r_1 + r_2 + r_3 + 10 {r_1 \choose 2} {r_2 \choose 2} {r_3 \choose 2}$
such minimal primes.
Since $\ChJ{t}$ is a radical monomial ideal,
it has no embedded primes.

\section{sectjuliacheck}{Related ideals II:
generated by $g_{i,a,b}$ when $d(a,b) = 2,3$}

We still use $d_t(a,b) = \#\{i \in [t]: a_i \not = b_i\}$
and $D_t(a,b) = \{i \in [t]: a_i \not = b_i\}$.

\defin[defJcheck]
Let $t \le n$
and let $\CcJ{t}$ be the ideal generated by $g_{i,a,b}$
with $a, b \in N$ and $i \in D_t(a,b)$
such that either $d(a,b) = d_t(a,b) = 3$ or $d(a,b) = 2$.
\endb

\thm[thmminprimescheck]
Let ${\bf \check S}$ be the collection of
all $t$-signed subsets $S$ of $N$
(see \ref{defswitchable})
such that for any $a,b \in N$,
if $d(a,b) = d_t(a,b) = 3$,
then at least one of $a$ and $b$ is not in $S$.
Then the prime ideals that are minimal over $\CcJ{t}$
are the minimal ideals in $\{\CQ{t}{S}: S \in {\bf \check S}\}$
(see \ref{defadmPQ}).
\endb

\proof
Let $\CQ{t}{S}$ be in $\bf \check S$.
By \ref{propCJtinCQ},
$\CJ{t} \subseteq \CQ{t}{S}$,
and by \ref{lemmamonogens} and by the $d_t, d$-conditions on $S$,
$\ChJ{t} \subseteq \CQ{t}{S}$.
Thus $\CcJ{t} \subseteq \CQ{t}{S}$.
By \ref{thmGBprime},
$\CQ{t}{S}$ is a prime ideal.

Note that $\CcJ{t} = \CJ{t} + \ChJ{t}$.
Let $P$ be a prime ideal minimal over $\CcJ{t}$.
Let $S$ be the set of all $a \in N$ such that $x_a \not \in P$.
As in the proof of \ref{thmCJminprimesone},
$S$ is $t$-switchable.
Thus by \ref{thmGBprime},
$\CQ{t}{S}$ is a prime ideal,
and by \ref{propCJtinCQ},
$\CQ{t}{S}$ contains $\CJ{t}$.
Since $P$ contains $\ChJ{t}$,
then at least one of $a$ and $b$ is not in $S$
whenever $d(a,b) = d_t(a,b) = 3$.
Again as in the proof of \ref{thmCJminprimes},
$S$ is $t$-signed,
and $\CcJ t \subseteq \CQ t S \subseteq P$.
Thus $S \in {\bf \check S}$,
and $\CQ t S = P$.
If $\CQ t T \subseteq \CQ t S$ for some $T \in {\bf \check S}$,
then by the first paragraph and minimality of $P$,
$\CQ t T = \CQ t S = P$ is minimal in ${\bf \check S}$.

Now assume that $\CQ{t}{S}$ is minimal in $\bf \check S$.
By the first paragraph,
$\CcJ{t} \subseteq \CQ{t}{S}$,
and $\CQ{t}{S}$ is a prime ideal.
Let $P$ be a prime ideal that is minimal over $\CcJ{t}$
and is contained in $\CQ{t}{S}$.
By the already established part,
$P = \CQ{t}{T}$ for some $T \in {\bf \check S}$.
Since $\CQ{t}{T} \subseteq \CQ{t}{S}$,
then by the minimality,
$\CQ{t}{T} = \CQ{t}{S}$.
\qed

\remark
We can give a more precise combinatorial description of the sets $S$
in the theorem above.
Let $S_0$ be a connected component of $S$.
Then $S_0$ is also in $\bf \check S$.
Suppose that for all $a, b \in S_0$,
$d_t(a,b) \le 2$.
We claim that $S_0$ is contained in a subset of $N$
in which all except some two coordinates in $[t]$ are fixed.
If not,
there exist $a, b, c \in S_0$
such that $a$ and $b$ differ in coordinates $i$ and $j$ in $[t]$,
$a$ and $c$ differ in coordinates $k$ and $l$ in $[t]$,
and $\{i,j,k,l\}$ has at least three elements.
Since $d_t(b,c) \le 2$,
necessarily $\#\{i,j,k,l\} = 3$.
Say $i = l$,
so that $b$ and $c$ differ in positions $j$ and $k$.
But then $b' = \switch(i,b,a)$ differs from $c$ in positions $i, j$ and $k$.
By~\ref{rmktadmswitcht},
$b', c' = \switch(\{i,j,k\},c,b') \in S_0$,
but $d_t(c,c') = d(c,c') = 3$,
which contradicts
that $S_0 \in {\bf \check S}$.
Thus indeed $S_0$ is contained in a subset of $N$
in which all except some two coordinates in $[t]$ are fixed.

\example
Let $t = n = 3$ and let $S \in {\bf \check S}$
be such that $\CQ{t}{S}$ is minimal over $\CcJ{t}$.
Since $t = n = 3$,
$d_t(a,b) = d(a,b) \le 3$ for all $a, b \in S$.
But since $S \in {\bf \check S}$,
$d(a,b) \not = 3$.
Thus for all $a, b \in S$,
$d(a,b) \le 2$.
Since $t = 3$,
by switchability,
all elements of $S$ are connected.
Thus by the remark,
there exists $k$ such that all elements of $S$ have the same $k$th coordinate.
Let $\{1,2,3\} = \{i,j,k\}$.
\ritem
If also the $j$th coordinate of all elements of $S$ is fixed,
then the minimal $\CQ t S$ can only be achieved
if $S$ is the whole line segment with those fixed $j$th and $k$th coordinates.
There are $r_1 \cdot r_2 + r_1 \cdot r_3 + r_2 \cdot r_3$ such sets.
\ritem
Now suppose that $S$ is not contained in a line segment.
Since $S$ is $t$-signed,
by \ref{lmtadm2},
$S$ is contained in a $2 \times 2$ submatrix
parallel to a coordinate plane.
To achieve a minimal $\CQ t S$,
$S$ then consists of all four points of that submatrix.
There are $r_1 {r_2 \choose 2} {r_3 \choose 2}
+ r_2 {r_1 \choose 2} {r_3 \choose 2}
+ r_3 {r_1 \choose 2} {r_2 \choose 2}$ such primes.
\ritem
Suppose that $r_i = 2$.
Then a prime ideal $\CQ t S$ of the type as in part (1)
strictly contains a prime ideal of part (2),
so not all $S$ in parts (1) and (2) correspond to minimal primes.

Note that there are no further containments among these prime ideals,
so that in total there are
$(r_1 \cdot r_2) \delta_{r_3 > 2}
+ (r_1 \cdot r_3) \delta_{r_2 > 2}
+ (r_2 \cdot r_3) \delta_{r_1 > 2}
+ r_1 {r_2 \choose 2} {r_3 \choose 2}
+ r_2 {r_1 \choose 2} {r_3 \choose 2}
+ r_3 {r_1 \choose 2} {r_2 \choose 2}$ minimal primes,
where $\delta_C$ is $1$ if condition $C$ is true and is $0$ otherwise.

\references

\bitem{AR}
N.~Ay and J.~Rauh, Robustness and conditional independence ideals,
(2011) {\tt arXiv:1110.1338}.

\bitem{BV}
W. Bruns and U. Vetter,
{\bkt Determinantal Rings},
Lecture Notes in Mathematics no. 1327,
Springer-Verlag, 1988.

\bitem{ESbinom}
D. Eisenbud and B. Sturmfels,
Binomial ideals,
{\it Duke Math. J.} {\bf 84} (1996), 1--45.

\bitem{Fink}
A. Fink,
The binomial ideal of the intersection axiom for conditional probabilities,
{\it J. Algebraic Combin.} Online: {\bf DOI: 10.1007/s10801-010-0253-5}.

\bitem{GS}
D. Grayson and M. Stillman,
Macaulay2, a software system for research in algebraic geometry,
available at http://www.math.uiuc.edu/Macaulay2/.

\bitem{HHHKR}
J. Herzog, T. Hibi, F. Hreinsdottir, T. Kahle, and J. Rauh,
Binomial edge ideals and conditional independence statements,
{\it Adv. in Appl. Math.} {\bf 45} (2010) 317--333.

\bitem{Ka}
T.~Kahle,
Decompositions of binomial ideals in Macaulay2,
{\tt arXiv:1106.5968}.

\bitem{Ki}
G. A. Kirkup,
Minimal primes over permanental ideals,
{\it Trans. Amer. Math. Soc.} {\bf 360} (2008), 3751--3770.

\bitem{LS}
R. C. Laubenbacher and I. Swanson,
Permanental ideals,
{\it J. Symbolic Comput.} {\bf 30} (2000), 195--205.

\bitem{O}
M. Ohtani,
Graphs and ideals generated by some $2$-minors,
{\it Comm. Alg.} {\bf 39} (2011), 905--917.

\bitem{ST}
I. Swanson and A. Taylor,
Minimal primes of ideals arising from conditional independence statements,
{\tt arXiv:math.AC/11075604}.

\bitem{V}
L. G. Valiant,
The complexity of computing the permanent,
{\it Theoret. Comp. Sci.} {\bf 8} (1979), 189--201.

\end

\section{sectnotneededI}{misc relations}

\remark[rmkfgrels]
The following are some simple relations
on the generalized permanents and determinants
when $l \not \in L$:
$$
\eqalignno{
g_{L \cup \{l\}, a, b}
&= f_{L,a,b} + g_{l, \switch(L,a,b), \switch(L,b,a)} \cr
&= g_{L,a,b} - f_{l, \switch(L,a,b), \switch(L,b,a)}, \cr
f_{L \cup \{l\}, a, b}
&= f_{L,a,b} + f_{l, \switch(L,a,b), \switch(L,b,a)} \cr
&= g_{L,a,b} - g_{l, \switch(L,a,b), \switch(L,b,a)}. \cr
}
$$

\section{sectGB}{Mm and mM}

For the purposes of discussion below we introduce the following notion:
we enumerate in the lexicographic order all possible elements of
$[r_{t+1}] \times \cdots \times [r_n]$.
Thun we can think of each element $a \in N$
as a $(t+1)$-tuple
in $[r_1] \times \cdots \times [r_t] \times [r_{t+1} \cdots r_n]$.

\defin[deftilde]
When we think of $a$ this way,
we denote $a$ as $\tilde a$.
\endb

In this section we also require that $t < n$
(recall that $t = n$ gives the same relevant ideals
as $t = n-1$).

\defin[defbigdiff]
Let $t < n$ and $a, b \in N$.
We define the {\bf smallest difference index} of $a$ and $b$,
$l(a,b)$,
to be the least integer such that $\tilde a_l \not = \tilde b_l$;
and if there is no such~$l$,
we declare $l(a,b) = t+2$.
With this,
we define the {\bf big difference, $\bf \bbD(a,b)$} of $a$ and $b$
as follows.
If $\tilde a_{t+1} \not = \tilde b_{t+1}$
or if $\tilde a_{t+1} = \tilde b_{t+1}$ and $\tilde a_l > \tilde b_l$,
then
$$
D(a,b) =
\#\{i \in [t] :
\sgn(a_i - b_i) = \sgn(\tilde b_{t+1} - \tilde a_{t+1} + 0.5)\} - 1,
$$
and otherwise set $D(a,b) = D(b,a)$.
\endb

Note that if $a = b$,
or more generally, if $l(a,b) \ge t+1$,
then $\bD(a,b) = 0$.

\lemma[lmDadditive]
Let $a, b \in S$ and $K \subseteq \{i \in [t] : a_i \not = b_i\}$.
Then $(\#K + \bD(a,b)) \cdot (\pl_S(a,b) - 1)$
has the same parity as
$\bD(\switch(K,a,b), \switch(K,b,a)) \cdot (\pl_S(a,b) - 1)$.
\endb

\proof
For brevity we denote $l = l(a,b)$,
$a_K = \switch(K,a,b)$,
and $b_K = \switch(K,b,a)$.
Note that $l(a,b) = l(a_K,b_K)$.

The conclusion is easy if $K = \emptyset$,
so we may assume that $K$ is not empty.
Thus $l \in [t]$.
Without loss of generality $a > b$ lexicographically,
i.e.,
$\tilde a_l > \tilde b_l$.

First note that
$$
\eqalignno{
&\hskip-4em\#\{i \in [t] :
\sgn((a_K)_i - (b_K)_i) = +\} \cr
&= \#\{i \in [t] :
{(a_K)}_i > {(b_K)}_i\} \cr
&= \#\{i \in K : a_i < b_i)\}
+ \#\{i \in [t] \setminus K : a_i > b_i)\} \cr
&= \#\{i \in K : a_i < b_i)\}
+ \#\{i \in K : a_i > b_i)\} \cr
&\hskip 2em
- 2 \#\{i \in K : a_i > b_i)\} \cr
&\hskip 2em
+ \#\{i \in K : a_i > b_i)\}
+ \#\{i \in [t] \setminus K : a_i > b_i)\} \cr
&= \#K - 2 \#\{i \in K : a_i > b_i)\} + \#\{i \in [t] : a_i > b_i)\} \cr
&= \#K - 2 \#\{i \in K : a_i > b_i)\}
+ \#\{i \in [t] : \sgn(a_i - b_i) = +\}, \cr
}
$$
so that
$\#\{i \in [t] : \sgn((a_K)_i - (b_K)_i) = +\}$
and $\#K + \#\{i \in [t] : \sgn(a_i - b_i) = +\}$
have the same parity.
Similarly
$\#\{i \in [t] : \sgn((a_K)_i - (b_K)_i) = -\}$
and $\#K + \#\{i \in [t] : \sgn(a_i - b_i) = -\}$
have the same parity.

In particular,
if $\tilde a_{t+1} \not = \tilde b_{t+1}$,
then
$$
\bD(a_K,b_K)
= \#\{i \in [t] :
\sgn((a_K)_i - (b_K)_i) =
\sgn(\widetilde {(b_K)}_{t+1} - \widetilde {(a_K)}_{t+1} + 0.5)\} - 1
$$
has the same parity as
$\#K + \#\{i \in [t] : \sgn(a_i - b_i) =
\sgn(\widetilde b_{t+1} - \widetilde a_{t+1} + 0.5)\} - 1
= \#K + \bD(a,b)$.

Now suppose that $\tilde a_{t+1} = \tilde b_{t+1}$.
Without loss of generality $a_l > b_l$.
If $l \not \in K$,
then $(a_K)_l > (b_K)_l$
and
$$
\bD(a_K,b_K) = \#\{i \in [t] : \sgn((a_K)_i - (b_K)_i) = +\} - 1
$$
has the same parity as
$\#K + \#\{i \in [t] : \sgn(a_i - b_i) = +\} - 1 = \#K + \bD(a,b)$.
Finally assume that $l \in K$.
Since the conclusion clearly holds if $\pl_S(a,b)$ is odd,
we now assume that $\pl_S(a,b)$ is even.
In this case $(a_K)_l < (b_K)_l$,
so that
$$
\bD(a_K,b_K) = \#\{i \in [t] : \sgn((b_K)_i - (a_K)_i) = +\} - 1
$$
has the same parity as
$$
\eqalignno{
\#K &+ \#\{i \in [t] : \sgn(b_i - a_i) = +\} - 1 \cr
&= \#K + \#\{i \in [t] : \sgn(a_i - b_i) = -\} - 1 \cr
&= \#K + \#\{i \in [t] : a_i \not = b_i\}
- \#\{i \in [t] : \sgn(a_i - b_i) = +\} + 1 - 2 \cr
&= \#K + d(a,b) - \bD(a,b) - 2 \cr
&= \#K + \pl_S(a,b) + \bD(a,b) - 2 \bD(a,b) - 2, \cr
}
$$
and since $\pl_S(a,b)$ is even,
$\bD(a_K,b_K)$ has the same parity as $\#K + \bD(a,b)$.
\qed

\lemma[lmredquad]
Let $S$ be a $t$-signed set,
and let $a, b \in S$ be connected.
Let $A \in N$ be such that
the first $l$ entries of $\tilde A$ are the maxima
of the corresponding entries of $\tilde a$ and $\tilde b$,
and the remaining entries of $\tilde A$ are the minima
of the corresponding entries of $\tilde a$ and $\tilde b$.
Let $B \in N$ be such that for all $i = 1, \ldots, t+1$,
$\{\tilde A_i, \tilde B_i\} = \{\tilde a_i, \tilde b_i\}$.
Then $A, B \in S$ are connected to $a, b$ in $S$,
and $x_a x_b$ reduces in the Gr\"obner basis sense
with respect to $\CG{t}{S}$ to
$$
(-1)^{\bD(a,b) \cdot (\pl_S(a,b)-1)} x_A x_B,
$$
and $x_A x_B$ is reduced with respect to $\CG{t}{S}$.
\endb

\proof
If $l(a,b) \ge t+1$,
then $x_a x_b$ is reduced,
$\bD(a,b) = 0$,
and $\{A, B\} = \{a,b\}$.

So we may assume that $l(a,b) \in [t]$.
Without loss of generality $\tilde a_{l(a,b)} > \tilde b_{l(a,b)}$.

If $\tilde a_{t+1} \le \tilde b_{t+1}$,
set $K = \{i > l(a,b): \tilde a_i > \tilde b_i\}$
$= \{i \in \{l(a,b), \ldots, t\}: a_i > b_i\}$.
Then $A = \switch(K,a,b)$, $B = \switch(K,b,a)$.
They are connected to $a$ by \ref{rmktadmswitcht}.
Thus $h_{S,K,a,b} \in \CG{t}{S}$.
Note that $\#K = \bD(a,b)$.
If $\bD(a,b) = 0$,
then as in the previous case $x_a x_b$ is reduced.
So we may assume that $\bD(a,b) > 0$.
Then $x_a x_b = \lt(h_{S,K,a,b})$ and 
$$
x_a x_b - h_{S,K,a,b}
= (-1)^{\#K \cdot (\pl_S(a,b)-1)} x_{\switch(K,a,b)} x_{\switch(K,b,a)}
= (-1)^{\bD(a,b) \cdot (\pl_S(a,b)-1)} x_A x_B,
$$
which proves that $x_a x_b$ reduces to the desired monomial
with respect to $\CG{t}{S}$.
Now suppose that there exists $h \in \CG{t}{S}$
with $\lt(h) = x_A x_B$.
Necessarily $h = h_{S,L,A,B}$
for some non-empty $L \subseteq [t]$
with $x_A x_B > x_{\switch(L,A,B)} x_{\switch(L,B,A)}$.
But by the definition of $A, B$,
this is impossible.
Thus $x_A x_B$ is reduced.

Now suppose that $\tilde a_{t+1} > \tilde b_{t+1}$.
Then set $K = \{l(a,b)\} \cup \{i > l(a,b): \tilde a_i < \tilde b_i\}$
$K = \{l(a,b)\} \cup \{i \in \{l(a,b), \ldots, t\}: a_i < b_i\}$.
Then $A = \switch(K,b,a)$, $B = \switch(K,a,b)$.
As before,
$A$ and $B$ are connected to $a$,
$h_{S,K,a,b} \in \CG{t}{S}$,
and $x_a x_b = \lt(h_{S,K,a,b})$.
But
$\#K = 1 + \#\{i \in [t]: \tilde a_i < \tilde b_i\}
= 2 + \bD(a,b)$ has the same parity as $\bD(a,b)$,
and the rest of the proof is similar to the previous case.
\qed

\defin[defmMMm]
With notation as in the lemma above,
we call $A$ the {\listvarsax-min} of $a,b$,
and we denote it by $\listvarsm(a,b)$,
and we call $B$ the {\bf min-Max} of $a,b$,
and we denote it by $\bf mM(a,b)$.
\endb

Note that $Mm(a,b) = Mm(b,a)$ and $mM(a,b) = Mm(b,a)$.
However,
$Mm(Mm(a,b),c)$ need not equal $Mm(Mm(a,c),b)$:
if $a = (2,2,2)$, $b = (1,3,3)$, $c = (2,2,4)$,
then
$$
\eqalignno{
Mm(Mm(a,b),c) &= Mm((2,2,2),(2,2,4)) = (2,2,4), \cr
Mm(Mm(a,c),b) &= Mm((2,2,4),(1,3,3)) = (2,2,3). \cr
}
$$
Or if $a = (2,2,3,3,2)$, $b = (1,3,2,2,3)$, $c = (2,2,3,2,4)$,
then
$$
\eqalignno{
Mm(Mm(a,b),c) &= Mm(2,2,2,2,2),(2,2,3,2,4)) = (2,2,3,2,2), \cr
Mm(Mm(a,c),b) &= Mm((2,2,3,3,2),(1,3,2,2,3)) = (2,2,2,2,2). \cr
}
$$

\cor[corredquad]
Let $S$ be a $t$-signed set,
let $a, b \in S$ be connected,
and let $K \subseteq [t]$.
Then
$$
Mm(\switch(K,a,b),\switch(K,b,a)) = Mm(a,b)
\hbox{ and }
mM(\switch(K,a,b),\switch(K,b,a)) = mM(a,b),
$$
and both $x_a x_b$ and
$(-1)^{\#K \cdot (\pl_S(a,b) - 1)} x_{\switch(K,a,b)} x_{\switch(K,b,a))}$
reduce with respect to $\CG{t}{S}$ to
$$
(-1)^{\bD(a,b) \cdot (\pl_S(a,b)-1)} x_{Mm(a,b)} x_{mM(a,b)}.
$$
\endb

\proof
Certainly the $Mm$ and $mM$ of $a,b$
and of $\switch(K,a,b),\switch(K,b,a)$ agree.
By \ref{lmredquad},
$x_a x_b$ reduces to the indicated term,
and
$(-1)^{\#K \cdot (\pl_S(a,b) - 1)} x_{\switch(K,a,b)} x_{\switch(K,b,a))}$
reduces to
$$
(-1)^{\#K \cdot (\pl_S(a,b) - 1)}
(-1)^{\bD(\switch(K,a,b),\switch(K,b,a)) \cdot (\pl_S(a,b)-1)}
x_{Mm(a,b)} x_{mM(a,b)}.
$$
By \ref{lmDadditive},
$\#K \cdot (\pl_S(a,b) - 1)
+ \bD(\switch(K,a,b),\switch(K,b,a)) \cdot (\pl_S(a,b)-1)$
has the same parity as
$\bD(a,b) \cdot (\pl_S(a,b)-1)$,
which finishes the proof of the corollary.
\qed